\documentclass[journal,10pt,onecolumn]{IEEEtran}
\usepackage[noadjust]{cite}
\usepackage[pdftex]{graphicx}
\usepackage{amsmath,amsfonts,amssymb}
\usepackage{dsfont}
\usepackage{algorithmic}
\usepackage{array}
\usepackage[caption=false,font=footnotesize]{subfig}
\usepackage{stfloats}
\usepackage{bigints}
\usepackage{url}
\usepackage{color}
\newtheorem{proposition}{Proposition}
\newtheorem{lemma}{--- Lemma}
 
\begin{document}

\title{Riemannian Gaussian Distributions on the \\ Space of Symmetric Positive Definite Matrices}

\author{Salem~Said, Lionel~Bombrun, Yannick Berthoumieu, \& Jonathan H. Manton,~\IEEEmembership{Senior~Member,~IEEE}}

\maketitle

\begin{abstract}
\, Data which lie in the space $\mathcal{P}_{m\,}$, of $m \times m$ symmetric positive definite matrices, (sometimes called \textit{tensor data}), play a fundamental role in applications including medical imaging, computer vision, and radar signal processing. An open challenge, for these applications, is to find a class of probability distributions, which is able to capture the statistical properties of data in $\mathcal{P}_{m\,}$, as they arise in real-world situations. The present paper meets this challenge by introducing \textit{Riemannian Gaussian distributions} on $\mathcal{P}_{m\,}$. Distributions of this kind were first considered by Pennec in $2006$. However, the present paper gives an exact expression of their probability density function for the first time in existing literature. This leads to two original contributions. First, a detailed study of statistical inference for Riemannian Gaussian distributions, uncovering the connection between maximum likelihood estimation and the concept of Riemannian centre of mass, widely used in applications. Second, the derivation and implementation of an expectation-maximisation algorithm, for the estimation of mixtures of Riemannian Gaussian distributions. The paper applies this new algorithm, to the classification of data in $\mathcal{P}_{m\,}$, (concretely, to the problem of texture classification, in computer vision), showing that it yields significantly better performance, in comparison to recent approaches. 
\end{abstract}
\begin{IEEEkeywords}
Symmetric positive definite matrices, tensor, Riemannian metric, Gaussian distribution, expectation-maximisation, texture
\end{IEEEkeywords}
\section{Introduction}
It has been known for quite some time, in fields ranging from multivariate statistics and information geometry~\cite{statistics}\cite{atkinson}, to matrix analysis~\cite{bhatia}, harmonic analysis and number theory~\cite{helgason}\cite{terras2}, that the space $\mathcal{P}_{m\,}$, of $m \times m$ symmetric positive definite matrices, can be equipped with a Riemannian metric, which gives it the structure of a Riemannian homogeneous space of negative curvature. In the present paper, this Riemannian metric is called the Rao-Fisher metric, (another popular name is affine-invariant metric), and is the subject of Section \ref{sec:geo}. 

During the past ten years, in response to the need for effective methods of processing data which lie in the space $\mathcal{P}_{m\,}$, the Rao-Fisher metric, (in addition to an alternative so-called log-Euclidean metric~\cite{pennec3}), has been widely used in engineering applications, which include medical imaging~\cite{pennec2}\cite{congedo}, continuum mechanics~\cite{moakher2}, radar signal processing~\cite{arnaudon}\cite{arnaudon1}, and
computer vision~\cite{image1,image2,dong,tuzel,caseiro}. These applications have evolved useful techniques for analysis, representation and classification of data which lie in the space $\mathcal{P}_{m\,}$. An introduction to these techniques may be found in~\cite{moakher_rev}. 

However, the literature still lacks a probabilistic model, rigorously defined yet tractable, which is able to represent the statistical variability of data in $\mathcal{P}_{m\,}$, and to bring the machinery of statistical inference, (estimation and hypothesis testing), to bear on such data. The present paper develops such a probabilistic model, by introducing a new class of probability distributions on the space $\mathcal{P}_{m\,}$, called \textit{Riemannian Gaussian distributions}. These distributions are the subject of Section \ref{sec:gaussian}, below.

A Riemannian Gaussian distribution, denoted $G(\bar{Y},\sigma)$, depends on two parameters, $\bar{Y} \in \mathcal{P}_m$ and $\sigma > 0$. The expression of its probability density function, generalising that of a Gaussian distribution on the Euclidean space $\mathbb{R}^{\scriptscriptstyle p}$, is given by
\begin{equation} \label{eq:introduction2}
p(Y | \,\bar{Y},\sigma) = \frac{1}{\zeta(\sigma)} \, \exp \left[ - \frac{d^{\,2}(Y,\bar{Y})}{2\sigma^2}\,\right]
\end{equation} 
Here, $d : \mathcal{P}_m \times \mathcal{P}_m \rightarrow \mathbb{R}_+$ is Rao's Riemannian distance, and the density is with respect to the Riemannian volume element of $\mathcal{P}_{m\,}$, which is henceforth denoted $dv(Y)$, (for required mathematical definitions, see Paragraph \ref{subsec:metric}). In comparison to a Gaussian distribution on $\mathbb{R}$, the normalising factor $\zeta(\sigma)$ plays the role of the factor $\sqrt{\,2\pi\,\sigma^{\scriptscriptstyle 2}}\,$.

It will be seen that Riemannian Gaussian distributions provide a statistical foundation for the concept of Riemannian centre of mass, (also called geometric mean, Riemannian mean, or Fr\'echet mean), currently essential to many applications~\cite{afsari}\cite{moakher1}. Recall the Riemannian centre of mass
$\hat{Y}_{\scriptscriptstyle N}$, of data points $Y_1, \ldots, Y_{\scriptscriptstyle N}$ in $\mathcal{P}_{m\,}$, is the unique global
minimiser of $\mathcal{E}_{\scriptscriptstyle N} : \mathcal{P}_m \rightarrow \mathbb{R}_+$,  
\begin{equation} \label{eq:introduction1}
 \mathcal{E}_{\scriptscriptstyle N}(Y) = \frac{1}{\strut N}\sum^{\scriptscriptstyle N}_{n=1} d^{\,2}(Y,Y_n)
\end{equation}
where, again, $d : \mathcal{P}_m \times \mathcal{P}_m \rightarrow \mathbb{R}_+$  is Rao's Riemannian distance. Since $\hat{Y}_{\scriptscriptstyle N}$ minimises the sum of squares of distances to the points $Y_1, \ldots, Y_{\scriptscriptstyle N}$, it is widely viewed as a representative, average, or mode of these points. 

Distributions of the form (\ref{eq:introduction2}) were considered by Pennec, who defined them on general Riemannian manifolds~\cite{pennec1}. However, in existing literature, their treatment remains incomplete, as it is based on asymptotic formulae, valid only in the limit where the parameter $\sigma$ is small, see~\cite{pennec1,lenglet,vemuri}. In addition to being inexact, such formulae are quite difficult, both to evaluate and to apply. These issues, (lack of an exact expression and difficulty of application), are overcome in the following. 

Indeed, Section \ref{sec:gaussian} opens with Proposition \ref{prop:z}, which gives an exact expression of the normalising factor $\zeta(\sigma)$, appearing in (\ref{eq:introduction2}). This confirms the important property that this factor does not depend on the parameter $\bar{Y}$. Proposition \ref{prop:z} is followed by Propositions \ref{prop:transf} and \ref{prop:gaussianpolar}, which together define a general method for sampling from a given Riemannian Gaussian distribution. These three propositions fully characterise Riemannian Gaussian distributions both theoretically and in view of computer simulations. 
  
Propositions \ref{prop:mle} and \ref{prop:ybar} describe the relationship between Riemannian Gaussian distributions and the concept of Riemannian centre of mass. Proposition \ref{prop:mle} states that the maximum likelihood estimate of the parameter $\bar{Y}$ of a distribution $G(\bar{Y},\sigma)$, based on samples $Y_1, \ldots, Y_{\scriptscriptstyle N}$ from this distribution, is equal to the Riemannian centre of mass $\hat{Y}_{\scriptscriptstyle N\,}$, global minimiser of (\ref{eq:introduction1}). In short, for Riemannian Gaussian distributions, \textit{maximum likelihood is equivalent to Riemannian centre of mass}. 

Proposition \ref{prop:ybar} is the asymptotic counterpart of Proposition \ref{prop:mle}, corresponding to the limit $N \rightarrow \infty\,$. It states that the parameter $\bar{Y}$ of a distribution $G(\bar{Y},\sigma)$ is the Riemannian centre of mass of this distribution, (see~\cite{afsari}, for the concept of Riemannian centre of mass of a probability distribution). This means that $\bar{Y}$ is the unique global minimiser of  $\mathcal{E} : \mathcal{P}_m \rightarrow \mathbb{R}_+$, 
\begin{equation} \label{eq:introduction3}
 \mathcal{E}(Y) = \int_{\mathcal{P}_m} d^{\,2}(Y,Z)\, p(Z|\, \bar{Y},\sigma) \, dv(Z)
\end{equation}
Note that, by the law of large numbers, $\mathcal{E}_{\scriptscriptstyle N}(Y) \rightarrow \mathcal{E}(Y)$ as $N \rightarrow \infty\,$, so (\ref{eq:introduction3}) is the asymptotic counterpart of (\ref{eq:introduction1}). Propositions \ref{prop:mle} and \ref{prop:ybar} allow for the concept of Riemannian centre of mass to be studied within the framework of parametric statistical inference. Compare to~\cite{bhatta1}\cite{bhatta2}, which have a non-parametric approach. 

This is the object of Proposition \ref{prop:hypothesis}, which deals with hypothesis testing, (testing against a given value of $\bar{Y}$), and of Propositions \ref{prop:consistency} and \ref{prop:normality}, which establish that $\hat{Y}_{\scriptscriptstyle N}$ is asymptotically normally distributed about $\bar{Y}$, in the limit $N \rightarrow \infty\,$, and give its asymptotic covariance matrix, (this makes it possible to construct asymptotic confidence regions for $\bar{Y}$). These results are intended to assist users wishing to assign a statistical significance to the Riemannian centre of mass $\hat{Y}_{\scriptscriptstyle N\,}$, computed from data points $Y_1, \ldots, Y_{\scriptscriptstyle N\,}$. 

Before going on, note that, in Section \ref{sec:gaussian}, $\hat{Y}_{\scriptscriptstyle N}$ is called by the different name of empirical Riemannian centre of mass, in order to distinguish it from $\bar{Y}$. Precisely, $\,\hat{Y}_{\scriptscriptstyle N}$ is an estimate, while $\bar{Y}$ is a parameter. 

The use of Riemannian Gaussian distributions in the representation and classification of data in $\mathcal{P}_{m\,}$ is considered in Section \ref{sec:classif}. 
This section is motivated by the idea that the class of \textit{mixtures of Riemannian Gaussian distributions} is expected to be sufficiently
rich, in order to represent the statistical distribution of data in $\mathcal{P}_m$ which arise in real-world applications, (experimental verification of this idea is still ongoing~\cite{gretsi}). A mixture of Riemannian Gaussian distributions is a probability distribution on $\mathcal{P}_m$ which has probability density function, 
\begin{equation} \label{eq:introduction4}
   p(Y) = \sum^{\scriptscriptstyle M}_{\mu \,= 1} \varpi_\mu \times p(Y | \,\bar{Y}_\mu,\sigma_\mu) 
\end{equation}
where $\varpi_1, \ldots, \varpi_{\scriptscriptstyle M}$ are non-zero positive weights which satisfy  $\varpi_1 + \ldots + \varpi_{\scriptscriptstyle M} = 1$, and where each density $p(Y | \,\bar{Y}_\mu,\sigma_\mu)$ is given by (\ref{eq:introduction2}), with $\bar{Y}_\mu \in \mathcal{P}_m$ and $\sigma_\mu > 0$.  

The estimation of mixtures of Riemannian Gaussian distributions is treated in Paragraph \ref{subsec:em}, which derives a new EM (expectation-maximisation) algorithm for computing maximum likelihood estimates of mixture parameters $\vartheta = \lbrace (\varpi_\mu,\bar{Y}_\mu,\sigma_\mu)\,; $
 $\mu = 1, \ldots, M \rbrace$. This algorithm appears as a generalisation, to the context of the Riemannian geometry of the space
$\mathcal{P}_{m\,}$, of EM algorithms currently used for estimation of mixtures of parametrised probability distributions on a Euclidean space, see~\cite{mixtures}\cite{mixture}, for example.

The classification of data in $\mathcal{P}_{m\,}$, using mixtures of Riemannian Gaussian distributions, is developed in Paragraph \ref{subsec:laplace}. The proposed approach aims to improve upon recent methods, used in~\cite{congedo}\cite{arnaudon}. While these are based on a purely geometric nearest-neighbor classification rule, Paragraph \ref{subsec:laplace} proposes a new statistical \textit{Bayes optimal classification rule}, which follows from the mixture distribution representation (\ref{eq:introduction4}).

The paper closes with Section \ref{sec:numerical}, which presents a numerical experiment showing the improvement in the rate of successful classification, obtained by the new classification rule of Paragraph \ref{subsec:laplace}, over the nearest-neighbour rule of~\cite{congedo}\cite{arnaudon},
and in comparison to a benchmark classification rule using a Wishart classifier~\cite{nielsen}\cite{wishartclass}. This numerical experiment is carried out on real data from the Vision Texture image database (VisTex)~\cite{vistex}.

Throughout the following, it should be kept in mind that the results of Sections \ref{sec:gaussian} and \ref{sec:classif} rely in a crucial way on the fact that the space $\mathcal{P}_{m\,}$, equipped with the Rao-Fisher metric, is a Riemannian homogeneous space of negative curvature. In particular, these results have no meaningful counterpart using alternative Riemannian metrics, such as the above-mentioned log-Euclidean metric. 

On the other hand, they generalise immediately to any Riemanniann homogeneous space of negative curvature. Spaces of this kind include the space of univariate normal distributions~\cite{entropy}, the space of Toeplitz autocovariance matrices~\cite{arnaudon}, and the space of Block-Toeplitz autocovariance matrices which have Toeplitz blocks~\cite{jeuris1}\cite{jeuris2}, with further examples given in~\cite{barbaresco1}\cite{barbaresco}. The definition and properties of Riemannian Gaussian distributions can be generalised to any of these spaces. This generalisation is the focus of ongoing work.

To place the present paper in this general context, consider briefly the definition of Riemannian Gaussian distributions on an arbitrary Riemannian homogeneous space of negative curvature, $\mathcal{M}$. This gives an abstract picture of Riemannian Gaussian distributions, which will be made concrete as of Section \ref{sec:geo}, by putting $\mathcal{M} = \mathcal{P}_{m\,}$. The definition is based on three essential ingredients\,: the existence of an invariant distance, the existence of invariant integrals, and the existence and uniqueness of Riemannian centres of mass. In the special case $\mathcal{M} = \mathcal{P}_{m\,}$, these will be formulated precisely in Propositions \ref{prop:metric}, \ref{prop:integral} and \ref{prop:curvature} of Section \ref{sec:geo}.
 
The space $\mathcal{M}$ is equipped with a Riemannian metric, which defines a Riemannian distance $d : \mathcal{M} \times \mathcal{M} \rightarrow \mathbb{R}_+$. The distance between two points $x,y  \in \mathcal{M}$ is then $d(x,y)$. To say that $\mathcal{M}$ is a Riemannian homogeneous space means two things \cite{helgason}\cite{terras2}.

First, a group of transformations $G$ acts transitively on $\mathcal{M}$. That is, each element $g \in G$ defines a transformation of $\mathcal{M}$, mapping each point $x \in \mathcal{M}$ to its image $x \cdot g$. Moreover, for any two points $x, y \in \mathcal{M}$ there exists $g \in G$ such that $y = x\cdot g$. Second, the Riemannian distance is invariant under the group of transformations $G$, 
$$
\mbox{invariant distance : }\hspace{0.5cm}  d(x,y) = d(x \cdot g,y\cdot g) \hspace{0.75cm} x,y \in \mathcal{M} \, , \, g \in G
$$
The notation $x \cdot g$ stands for right action of $G$ on $\mathcal{M}$. That is, $x \cdot (g_{\scriptscriptstyle 1}g_{\scriptscriptstyle 2}) = (x \cdot g_{\scriptscriptstyle 1})\cdot g_{\scriptscriptstyle 2\,}$, for $g_{\scriptscriptstyle 1\,},g_{\scriptscriptstyle 2\,} \in G$. The choice of right action over left action, (which would be denoted $g\cdot x$, instead of $x\cdot g$), is a matter of convention, which is here adopted throughout the paper.

A corollary of the invariance of Riemannian distance is that it is possible to define invariant integrals. Precisely, the Riemannian metric of $\mathcal{M}$ defines a Riemannian volume element, $dv(x)\;$\cite{terras2}. Integrals with respect to this volume element are invariant under the group of transformations $G$,
$$
\mbox{invariant integrals : }\hspace{0.5cm}  \int_{\scriptscriptstyle \mathcal{M}}\, f(x)dv(x) = \int_{\scriptscriptstyle \mathcal{M}}\, f(x\cdot g)dv(x) \hspace{0.75cm} g \in G
$$
This property means that the integral of any function $f : \mathcal{M} \rightarrow \mathbb{R}$ is equal to the integral of this function composed with a transformation $x \mapsto x\cdot g\,$. Intuitively, it is a direct result of the invariance of distance, through the basic relationship between the concepts of length and volume. 

In the special case $\mathcal{M} = \mathcal{P}_{m\,}$, the existence of invariant distance and of invariant integrals are formulated in equations (\ref{eq:invmetric2}) and (\ref{eq:integration1}) below. In all generality, these two properties supply the expression of the probability density function of a Riemannian Gaussian distribution on $\mathcal{M}$,
$$
\mbox{Riemannian Gaussian p.d.f. : }\hspace{0.5cm} p(x|\,\bar{x},\sigma) = \frac{1}{\zeta(\sigma)} \, \exp \left[ - \frac{d^{\,2}(x,\bar{x})}{2\sigma^2}\,\right]
$$
where $\bar{x} \in \mathcal{M}$ and $\sigma > 0$ are parameters, and the density is with respect to the Riemannian volume element $dv(x)$.

 The crucial feature of this expression of the probability density is that the normalising factor $\zeta(\sigma)$ does not depend on the parameter $\bar{x}$, but only on the parameter $\sigma$. This follows from the existence of invariant distance and of invariant integrals, as can be seen by repeating word for word the proof of item (i) of Proposition \ref{prop:z}, Section \ref{sec:gaussian}.

 This feature implies the special relationship between Riemannian Gaussian distributions and Riemannian centres of mass. Namely, for these distributions, \textit{maximum likelihood is equivalent to Riemannian centre of mass}. Indeed, if $x_{\scriptscriptstyle 1}, \ldots, x_{\scriptscriptstyle N}$ are samples from the above probability density, the resulting log-likelihood function is
$$
\mbox{log-likelihood : }\hspace{0.5cm} \sum^{\scriptscriptstyle N}_{\scriptscriptstyle n=1} \, \log \, p(x_n|\,\bar{x},\sigma) = \,- N \log \zeta(\sigma) \, - \, \frac{1}{\mathstrut \sigma^{\scriptscriptstyle 2}} \, \sum^{\scriptscriptstyle N}_{\scriptscriptstyle n=1} \, d^{\,2}(x_n,\bar{x})
$$
Since the first term on the right hand side does not depend on $\bar{x}$, it follows that the maximum likelihood estimate $\hat{x}_{\scriptscriptstyle N}$ of the parameter $\bar{x}$ should minimise the sum of squares of distances to the samples $x_{n\,}$. In other words, it should be a Riemannian centre of mass of these samples, (compare to formula (\ref{eq:introduction1}), above). In the special case $\mathcal{M} = \mathcal{P}_{m\,}$, this reasoning will appear in the proof of Proposition \ref{prop:mle}, Section \ref{sec:gaussian}.

At this point, the fact that $\mathcal{M}$ has negative curvature becomes important. It is this fact that guarantees the existence and uniqueness of Riemannian centres of mass in $\mathcal{M}$~\cite{afsari}\cite{moakher1}, so that the maximum likelihood estimate $\hat{x}_{\scriptscriptstyle N}$ of the parameter $\bar{x}$ is well-defined and unique. In the special case $\mathcal{M} = \mathcal{P}_{m\,}$, this is stated in Proposition \ref{prop:curvature} of Section \ref{sec:geo}.

With regard to the maximum likelihood estimate $\hat{\sigma}_{\scriptscriptstyle N}$ of the parameter $\sigma$, a more detailed knowledge of the space $\mathcal{M}$ is required. Indeed, to establish existence and uniqueness of $\hat{\sigma}_{\scriptscriptstyle N\,}$, it is necessary to show that $\log \zeta(\sigma)$ is a strictly convex function of the ``natural parameter'' $\eta = -1/\sigma^{2\,}$. In the special case $\mathcal{M} = \mathcal{P}_{m\,}$, this follows from item (ii) of Proposition \ref{prop:z}, as stated in the proof of Proposition \ref{prop:mle}.

To obtain the same result for a general Riemannian homogeneous space of negative curvature, $\mathcal{M}$, it is necessary to decompose $\mathcal{M}$ into a direct product of symmetric spaces of Euclidean or of non-compact type, and consequently use integration formulas on these Symmetric spaces, which may be found in~\cite{helgason}. In a future submission, this will be carried out for all above mentioned examples, (\textit{e.g.}, spaces of autocovariance matrices with Toeplitz structure~\cite{arnaudon}, or Block-Toeplitz structure with Toeplitz blocks~\cite{jeuris1}\cite{jeuris2}).  However, the present paper will exclusively focus on the special case $\mathcal{M} = \mathcal{P}_{m\,}$, since it is relevant to a wider range of applications, and requires less mathematical background.

\section{Riemannian geometry of covariance matrices} \label{sec:geo}
Let $\mathcal{P}_m$ denote the space of all $m \times m$ real matrices $Y$ which are symmetric and strictly positive definite,
\begin{equation} \label{eq:pn}
   Y^\dagger - Y = 0 \hspace{1cm} x^\dagger Y x > 0 \mbox{ for all } x \in \mathbb{R}^m
\end{equation}
where $\dagger$ denotes the transpose. The present section is concerned with the Riemannian geometry of the space $\mathcal{P}_{m\,}$, when this space is equipped with the Rao-Fisher metric~\cite{helgason}\cite{terras2}. 

Paragraph \ref{subsec:metric} gives the expressions of the Rao-Fisher metric, and of the related Riemannian distance (Rao's distance), volume element, and geodesic curves. 

Paragraph \ref{subsec:geo} states the fundamental geometric properties of the space $\mathcal{P}_{m\,}$, in view of the development of following sections. Essentially, these properties reflect the fact that $\mathcal{P}_{m\,}$, (equipped with the Rao-Fisher metric), is a Riemannian homogeneous space of negative curvature.
\subsection{The Rao-Fisher metric\,: distance, geodesics and volume} \label{subsec:metric}
A Riemannian metric on $ \mathcal{P}_m$ is a quadratic form $ds^2(Y)$ which measures the \textit{squared length of a small displacement} $dY$, separating two elements $Y \in \mathcal{P}_m$ and $Y + dY \in \mathcal{P}_{m\,}$. Here, $dY$ is a symmetric matrix, since $Y$ and $Y + dY$ are symmetric, by (\ref{eq:pn}). The Rao-Fisher metric is the following~\cite{atkinson}\cite{terras2},
\begin{equation} \label{eq:metric}
   ds^2(Y) = \mathrm{tr}\, [Y^{-1}dY]^2
\end{equation}
where $\mathrm{tr}$ denotes the trace. When there is no loss of clarity, this is written $ds^2(Y) = \Vert dY \Vert^{2\,}$. 

The Rao-Fisher metric, like any other Riemannian metric on $\mathcal{P}_{m\,}$, defines a Riemannian distance $d : \mathcal{P}_m \times \mathcal{P}_m \rightarrow \mathbb{R}_+$. This is called Rao's distance, and is defined as follows~\cite{helgason}\cite{terras2}. Let $Y, Z \in \mathcal{P}_m$ and $c : [0,1] \rightarrow \mathcal{P}_m$ be a differentiable curve with $c(0) = Y$ and $c(1) = Z$. The length $L(c)$ of $c$ is defined by
\begin{equation} \label{eq:length}
   L(c)  \,= \, \int^1_0 ds(c(t)) \,= \,\int^1_0 \Vert \dot{c}(t) \Vert \,dt
\end{equation}
where $\dot{c}(t) = \frac{dc}{dt}\,$. Rao's distance $d(Y,Z)$ is the infimum of $L(c)$ taken over all differentiable curves $c$ as above.

When equipped with the Rao-Fisher metric, the space $\mathcal{P}_m$ is a Riemannian manifold of negative sectional curvature~\cite{helgason}\cite{terras2}. One implication of this property, (since $\mathcal{P}_m$ is also complete and simply connected), is that the infimum of $L(c)$ is realised by a unique curve $\gamma\,$, known as the geodesic connecting $Y$ and $Z$. The equation of this curve is the following,
\begin{equation} \label{eq:geodesic1}
   \gamma(t) = Y^{1/2} \,( Y^{-1/2} Z Y^{-1/2} )^{t} \,Y^{1/2} \hspace{1cm} \mbox{ for } t \in [0,1]
\end{equation}
Given expression (\ref{eq:geodesic1}), it is possible to compute $L(\gamma)$ from (\ref{eq:length}). This is precisely Rao's distance $d(Y,Z)$. It turns out~\cite{terras2},
\begin{equation} \label{eq:distance}
   d^{\,2}(Y,Z) = \mathrm{tr}\, [ \log(Y^{-1/2}Z Y^{-1/2}) ]^2
\end{equation}
All matrix functions appearing in (\ref{eq:geodesic1}) and (\ref{eq:distance}), (square root, elevation to the power $t$, and logarithm), should be understood as symmetric matrix functions. For this, see~\cite{higham}.  

Since the Rao-Fisher metric gives a means of measuring length, it can also be used to measure volume. Roughly, this is based on the elementary fact that ``the volume of a cube is the product of the lengths of its sides." The Riemannian volume element associated to the Rao-Fisher metric is the following~\cite{terras2}, 
\begin{equation} \label{eq:dv}
  dv(Y) = \det(Y)^{-\frac{m+1}{2}}\prod_{i \leq j} dY_{ij}
\end{equation}
where indices denote matrix entries. 

To close this paragraph, consider the expressions of $ds^2(Y)$ and $dv(Y)$, given by (\ref{eq:metric}) and (\ref{eq:dv}), in terms of polar coordinates.  The use of polar coordinates refers to the parametrisation of $Y \in \mathcal{P}_m$ by its eigenvalues and eigenvectors. For $(r_1, \ldots, r_m) \in \mathbb{R}^m$ and $U \in \mathrm{O}(m)$, where $\mathrm{O}(m)$ is the group of $m \times m$ real orthogonal matrices, let $\mathrm{diag}(e^{r})$ be the diagonal matrix with main diagonal $(e^{r_1},\ldots,e^{r_m})$, and
\begin{equation} \label{eq:polar}
  Y(r,U) = U^\dagger \, \mathrm{diag}(e^{r}) \, U
\end{equation}
Then, it is clear that $Y(r,U)$ belongs to $\mathcal{P}_{m\,}$. Conversely, given any $Y \in \mathcal{P}_{m\,}$, by writing down the spectral decomposition of $Y$, it is possible to find $(r,U)$ such that $Y = Y(r,U)$. 

Expressed in polar coordinates, $ds^2(Y)$ and $dv(Y)$ appear as follows,
\begin{equation} \label{eq:polarmetric}
ds^2(Y) = \sum^m_{j=1} dr^2_j +  8 \sum_{i < j} \sinh^2\left(\frac{r_i - r_j}{2}\right) \theta^2_{ij} 
\end{equation}  
\begin{equation} \label{eq:polardv}
 dv(Y) = 8^{\frac{m(m-1)}{4}} \det(\theta)  \prod_{i < j } \sinh\left(\frac{|r_i - r_j|}{2}\right) \prod^m_{i=1} dr_i 
\end{equation}
where $\theta_{ij} = \sum_k U_{jk} \, dU_{ik\,}$, and the notation $\det(\theta)$ is used for the exterior product,
\begin{equation} \label{eq:dettheta}
  \det(\theta) = \bigwedge_{i < j} \, \theta_{ij}
\end{equation}
Under a slightly different form, these expressions are given in~\cite{terras2} (Exercise $25$, Page $24$).

When using expression (\ref{eq:polardv}) to compute integrals with respect to the volume element $dv(Y)$, one should be careful that the correspondence between $Y$ and $(r,U)$ is not unique. In fact, even when $Y$ has distinct eigenvalues, there are $m!\, 2^m$ ways of choosing $(r,U)$. The factor $m!$ corresponds to all possible reorderings of $r_1, \ldots, r_{m\,}$, and the factor $2^m$ corresponds to the orientation of the columns of $U$, (multiplication by $+1$ or $-1$). Accordingly, for any function $f : \mathcal{P}_m \rightarrow \mathbb{R}$,
\begin{equation} \label{eq:integrationew}
  \int_{\mathcal{P}_m} f(Y) \, dv(Y) =  (m!\, 2^m)^{\scriptscriptstyle -1} \, \times \, 8^{\frac{m(m-1)}{4}} \, \int_{\mathrm{O}(m)}\, \int_{\mathbb{R}^m}  f(Y(r,U)) \, \det(\theta) \,
\prod_{i < j }  \sinh\left(\frac{|r_i - r_j|}{2}\right) \prod^m_{i=1} dr_i    
\end{equation}
where division by $m!\, 2^m$ cancels out the ambiguity in choosing $(r,U)$. This formula is given in~\cite{terras2}, (Proposition $3$, Page $43$). 
\subsection{The Rao-Fisher metric\,: geometric properties} \label{subsec:geo}
The Rao-Fisher metric, when introduced on the space $\mathcal{P}_{m\,}$, leads to many important geometric properties. Among these properties, the ones which will be used in following sections are here recalled and briefly explained. 

First, consider the fact that this metric turns the space $\mathcal{P}_m$ into a \textit{Riemannian homogeneous space} under the action of the linear group $\mathrm{GL}(m)$. Recall $\mathrm{GL}(m)$ is the group of $m \times m$ real invertible matrices. This group acts on $\mathcal{P}_m$ by congruence transformations, which are defined as follows \cite{terras2},  
\begin{equation} \label{eq:congruence}
  (Y,A) \mapsto Y\cdot A \hspace{1cm} Y\cdot A = A^\dagger Y A 
\end{equation}
for $Y \in \mathcal{P}_m$ and $A \in \mathrm{GL}(m)$. The notation $Y\cdot A$ is well suited, in view of the right action property\,: $Y\cdot (A_1A_2) = (Y\cdot A_1)\cdot A_{2\,}$, for $A_1,A_2 \in \mathrm{GL}(m)$.

The term  ``homogeneous space'' means that for any $Y,Z \in \mathcal{P}_m$ there exists $A \in \mathrm{GL}(m)$ such that $Y\cdot A = Z$. One possible choice of $A$ is $A = Y^{-1/2}Z^{1/2}$. Moreover, the term ``Riemannian homogeneous space" means that the Rao-Fisher metric and Rao's distance  remain invariant under the action of $\mathrm{GL}(m)$ on $\mathcal{P}_{m\,}$. This is stated in Proposition \ref{prop:metric}.
\vspace{0.1cm}
\begin{proposition}[Riemannian Homogeneous space] \label{prop:metric}
  For $Y,Z \in \mathcal{P}_m$ and $A \in \mathrm{GL}(m)$, the following hold,
\begin{subequations} \label{eq:invmetric}
\begin{equation} \label{eq:invmetric1}
  ds^2(Y) = ds^2(Y \cdot A) \hspace{0.75cm} ds^2(Y) = ds^2(Y^{-1})
\end{equation} 
where $ds^2(Y)$ is the Rao-Fisher metric (\ref{eq:metric}).
\begin{equation} \label{eq:invmetric2}
  d(Y,Z) = d(Y \cdot A,Z\cdot A) \hspace{0.75cm} d(Y,Z) = d(Y^{-1},Z^{-1})
\end{equation} 
where $d: \mathcal{P}_m \times \mathcal{P}_m \rightarrow \mathbb{R}_+$ is Rao's distance (\ref{eq:distance}).
\end{subequations}
\end{proposition}
\vspace{0.1cm}
\textbf{Proof\,: } Identities (\ref{eq:invmetric}) are well-known in the literature~\cite{helgason}\cite{terras2}. They state that congruence transformations, as well as matrix inversion, are \textit{isometries} of the space $\mathcal{P}_{m\,}$, equipped with the Rao-Fisher metric. 

By general arguments from Riemannian geometry, (\ref{eq:invmetric2}) is a direct result of (\ref{eq:invmetric1}). Here, as an illustration, is the proof
of the first identity in (\ref{eq:invmetric1}). Let $W = Y\cdot A$, so $dW = dY \cdot A$. Replace in (\ref{eq:metric}),
$$
ds^2(W) = \mathrm{tr} \, [(A^\dagger Y A)^{-1} A^\dagger dY A]^2 = \mathrm{tr} \, [ Y^{-1}dY]^2 
$$
which is just $ds^2(Y)$. For the second identity, recall that if $W = Y^{-1}$ then $dW = -WdYW$, see~\cite{higham}.\hfill $\blacksquare$ \\[0.15cm]
\indent The following Proposition \ref{prop:integral} states that integrals with respect to the Riemannian volume element $dv(Y)$ are invariant under the action of $\mathrm{GL}(m)$ on $\mathcal{P}_{m\,}$.
\begin{proposition}[Invariant integrals]  \label{prop:integral}
  For any function $f : \mathcal{P}_m \rightarrow \mathbb{R}$, and any $A \in \mathrm{GL}(m)$,
\begin{subequations} \label{eq:integration}
  \begin{equation} \label{eq:integration1}
      \int_{\mathcal{P}_m} f(Y) \, dv(Y) = \int_{\mathcal{P}_m} f(Y\cdot A) \, dv(Y) 
  \end{equation}
  \begin{equation} \label{eq:integration2}
      \int_{\mathcal{P}_m} f(Y) \, dv(Y) = \int_{\mathcal{P}_m} f(Y^{-1}) \, dv(Y) 
  \end{equation}
whenever these integrals exist.
\end{subequations}
\end{proposition}
\textbf{Proof\,: } This proposition is a corollary of Proposition \ref{prop:metric}. Intuitively, it holds because $dv(Y)$ is the Riemannian volume element of $ds^2(Y)$, and $ds^2(Y)$ is invariant under congruence transformations and inversion, according to Proposition \ref{prop:metric}. A simple proof may be found in~\cite{terras2}.\hfill $\blacksquare$ \\[0.15cm]
\indent The third and last property of the Rao-Fisher metric, which will be needed in the following, refers to the \textit{existence and uniqueness of Riemannian centres of mass}. 

Let $\pi$ be a probability distribution on $\mathcal{P}_{m\,}$. The variance function of $\pi$ is $\mathcal{E}_\pi : \mathcal{P}_m \rightarrow \mathbb{R}_+$,
\begin{equation} \label{eq:variance}
  \mathcal{E}_\pi(Y) = \int_{\mathcal{P}_m} d^{\,2}(Y,Z)\, d\pi(Z)
\end{equation}
Following~\cite{afsari}, a point $\bar{Y}_\pi \in \mathcal{P}_{m\,}$, which is a global minimiser of $\mathcal{E}_\pi$, is called a Riemannian centre of mass of $\pi$. 

The following Proposition \ref{prop:curvature} is a result of the property that the space $\mathcal{P}_{m\,}$, equipped with the Rao-Fisher metric, is a Riemannian manifold of negative curvature. 
\begin{proposition}[Riemannian centre of mass]  \label{prop:curvature}
  If $\pi$ is a probability distribution on $\mathcal{P}_{m\,}$, then $\pi$ has a unique Riemannian centre of mass $\bar{Y}_{\pi\,}$. Moreover, $\bar{Y}_\pi$ is the unique stationary point of the variance function $\mathcal{E}_{\pi\,}$. 
\end{proposition}
\vspace{0.2cm}
\textbf{Proof\,: } The proposition is a corollary of a general theorem in~\cite{afsari} (Theorem $2.1.$, Page $659$). This theorem states the proposition will hold, as soon as it can be shown that $\mathcal{P}_{m\,}$, equipped with the Rao-Fisher metric, is a Riemannian manifold of negative sectional curvature. However, this is a well known fact, which may be found in~\cite{helgason}\cite{terras2}. Moreover, in~\cite{lenglet} (Theorem $2.2.2.$, Page $428$), an explicit formula is given, for the Riemannian curvature tensor of $\mathcal{P}_{m\,}$, and is used through a direct calculation, to find the sectional curvature of $\mathcal{P}_m$ and check that it is negative. \hfill $\blacksquare$ 
\section{Riemannian Gaussian distributions} \label{sec:gaussian}
The main theoretical contribution of the present paper is to give an exact formulation of Riemannian Gaussian distributions. A Riemannian Gaussian distribution $G(\bar{Y},\sigma)$ is a probability distribution on $\mathcal{P}_{m\,}$, given by the probability density function, with respect to the Riemannian volume element (\ref{eq:dv}),
\begin{equation} \label{eq:gaussianpdf}
  p(Y | \,\bar{Y},\sigma) = \frac{1}{\zeta(\sigma)} \, \exp \left[ - \frac{d^{\,2}(Y,\bar{Y})}{2\sigma^2}\,\right]
\end{equation}
where $\bar{Y} \in \mathcal{P}_m$ and $\sigma > 0$ are parameters, and $d(Y,\bar{Y})$ is Rao's distance given by (\ref{eq:distance}). For brevity, the term Gaussian distribution will be used, instead of Riemannian Gaussian distribution. The present section is organised as follows.

Paragraph \ref{subsec:gaussiandef} is concerned with the definition of Gaussian distributions. In this paragraph, Proposition \ref{prop:z} gives an exact expression of the normalising factor $\zeta(\sigma)$, and Propositions \ref{prop:transf} and \ref{prop:gaussianpolar} define a general method for sampling from a given Gaussian distribution.

Paragraph \ref{subsec:gaussianstat} studies statistical inference problems for Gaussian distributions. In this paragraph, Proposition \ref{prop:mle} deals with maximum likelihood estimation of the parameters $\bar{Y}$ and $\sigma$, while Proposition \ref{prop:hypothesis} deals with the problem of hypothesis testing, (precisely, testing against a given value of $\bar{Y}$). 

Paragraph \ref{subsec:gaussianasymp} recovers the asymptotic properties of the maximum likelihood estimate of $\bar{Y}$. These properties, of consistency and asymptotic normality, are given in Propositions \ref{prop:consistency} and \ref{prop:normality}. They essentially rely on Proposition \ref{prop:ybar}, which states that the parameter $\bar{Y}$ of a distribution $G(\bar{Y},\sigma)$ is the Riemannian centre of mass of this distribution.

The idea of considering probability densities on $\mathcal{P}_{m\,}$, which depend on Rao's distance as in (\ref{eq:gaussianpdf}), is due to Pennec~\cite{pennec1}. What has been lacking, in order to make this idea generally applicable, is a better understanding of the normalising factor $\zeta(\sigma)$. The desire to achieve such an understanding is the starting point of the present section. 

\subsection{Definition and basic properties} \label{subsec:gaussiandef}
To define a Gaussian distribution $G(\bar{Y},\sigma)$, by means of the probability density function (\ref{eq:gaussianpdf}), it is necessary to have an exact expression of the normalising factor $\zeta(\sigma)$. This is given by the following Proposition \ref{prop:z}. 

Item (i) of this proposition confirms an important property of $\zeta(\sigma)$. Namely, this normalising factor does not depend on the parameter $\bar{Y}$. On the other hand, item (ii) states the expression of $\zeta(\sigma)$, formula (\ref{eq:z}).

Formula (\ref{eq:z}) may look somewhat complicated. However, for values of $m$ up to $m = 50$, it has been easily evaluated using Monte Carlo integration. In particular, this has made it possible to build tables of $\zeta(\sigma)$ as a function of $\sigma$. For $m = 2$, formula (\ref{eq:z}) yields the analytic expression
\begin{equation} \label{eq:z2}
 \zeta(\sigma) = (2 \pi)^{3/2} \, \sigma^2 \times e^{\sigma^2/4} \times \mathrm{erf}(\sigma/2)
\end{equation}
where $\mathrm{erf}$ denotes the error function~\cite{lebedev}.

In order to state Proposition \ref{prop:z}, consider the following notation. For $\bar{Y} \in \mathcal{P}_m$ and $\sigma > 0$, let $f(Y|\, \bar{Y},\sigma)$ be given by,
\begin{equation} \label{eq:f}
f(Y|\, \bar{Y},\sigma) = \exp \left[ - \frac{d^{\,2}(Y,\bar{Y})}{2\sigma^2}\,\right]
\end{equation}
Also, let $\zeta(\bar{Y},\sigma)$ be the integral,
\begin{equation} \label{eq:prez}
\zeta(\bar{Y},\sigma) = \int_{\mathcal{P}_m} f(Y|\, \bar{Y},\sigma) \, dv(Y)
\end{equation}
where $dv(Y)$ is the Riemannian volume element (\ref{eq:dv}).  
\vspace{0.1cm}
\begin{proposition}[Normalising factor]  \label{prop:z}
The following hold, \\
\begin{subequations} \label{eq:superz} 
(i) For any $\bar{Y} \in \mathcal{P}_m$ and $\sigma > 0$,
\begin{equation} \label{eq:zinvariance}
  \zeta(\bar{Y},\sigma) = \zeta(I,\sigma)
\end{equation}
where $I \in \mathcal{P}_{m\,}$ is the $m \times m$ identity matrix. \\
(ii) Let $\zeta(\sigma) = \zeta(I,\sigma)$. Then, 
\begin{equation} \label{eq:z}
\zeta(\sigma) = (m!\, 2^m)^{\scriptscriptstyle -1} \, \times \omega_m \, \times \, 8^{\frac{m(m-1)}{4}} \,  \int_{\mathbb{R}^m} e^{- (r^{\,2}_1 + \ldots + r^{\,2}_m)/2\sigma^2}  \, \prod_{i < j } \sinh\left(|r_i - r_j|/2\right) \, \prod^m_{i =1 } dr_i
\end{equation}
where $\omega_m$ is given by  
\begin{equation} \label{eq:omegam}
\omega_m = \frac{2^m \pi^{m^2/2}}{\strut \Gamma_m (m/2)}
\end{equation}
\end{subequations}
with $\Gamma_m$ the multivariate Gamma function, given in~\cite{muirhead}.
\end{proposition}
\vspace{0.2cm}
\textbf{Proof\,: }  Item (i) follows from Propositions \ref{prop:metric} and  \ref{prop:integral}, in \ref{subsec:geo}. Let $\bar{Y} = I\cdot A\,$, in the notation of  (\ref{eq:congruence}). Replacing (\ref{eq:invmetric2}) in (\ref{eq:f}), it follows that 
$$
f(Y|\, \bar{Y},\sigma) = f(Y\cdot A^{-1}| \, I,\sigma)
$$
Then, applying (\ref{eq:integration1}) to (\ref{eq:prez}), it follows that $\zeta(\bar{Y},\sigma)$ is given by
$$
\int_{\mathcal{P}_m} f(Y\cdot A^{-1}| \, I,\sigma) \, dv(Y) = \int_{\mathcal{P}_m} f(Y| \, I,\sigma) \, dv(Y)
$$
which is just $\zeta(I,\sigma)$. This proves (\ref{eq:zinvariance}). A similar proof of this item may be found in~\cite{vemuri}, (Theorem 2.1., Page $597$).

Item (ii) follows by using (\ref{eq:integrationew}) to compute the integral $\zeta(\sigma) = \zeta(I,\sigma)$. Indeed, if $Y = Y(r,U)$ as in (\ref{eq:polar}), then it follows from (\ref{eq:distance}), see~\cite{lenglet} (Theorem $2.2.3.$, Page $430$),
\begin{equation} \label{eq:zproof1}
 f(Y| \, I,\sigma) = e^{- (r^{\,2}_1 + \ldots + r^{\,2}_m)/2\sigma^2}  
\end{equation}
Since this does not depend on $U$, replacing in (\ref{eq:integrationew}) gives
\begin{eqnarray}
\nonumber \zeta(\sigma) =  (m!\, 2^m)^{\scriptscriptstyle -1} \, \times \, 8^{\frac{m(m-1)}{4}} \, \int_{\mathrm{O}(m)}\, \int_{\mathbb{R}^m}  e^{- (r^{\,2}_1 + \ldots + r^{\,2}_m)/2\sigma^2} \, \det(\theta) \,
\prod_{i < j }  \sinh\left(\frac{|r_i - r_j|}{2}\right) \prod^m_{i=1} dr_i \\[0.1cm]
= (m!\, 2^m)^{\scriptscriptstyle -1} \, \times \, 8^{\frac{m(m-1)}{4}} \, \int_{\mathrm{O}(m)}\, \det(\theta) \, \int_{\mathbb{R}^m}  e^{- (r^{\,2}_1 + \ldots + r^{\,2}_m)/2\sigma^2} 
\prod_{i < j }  \sinh\left(\frac{|r_i - r_j|}{2}\right) \prod^m_{i=1} dr_i \,
\end{eqnarray}
Recall that, (see \cite{muirhead}, Page $71$), 
\begin{equation} \label{eq:zinproof3}
   \int_{\mathrm{O}(m)} \det(\theta)\, = \, \omega_m
\end{equation}
Now, (\ref{eq:z}) follows immediately. \hfill $\blacksquare$ \\[0.1cm]
\indent The following Propositions \ref{prop:transf} and \ref{prop:gaussianpolar} define a general method for sampling from a given Gaussian distribution, which will be described at the end of this section. To begin, Proposition \ref{prop:transf} gives the transformation properties of Gaussian distributions.

For the statement of this proposition, recall that if $X$ is a random variable and $\mu$ a probability distribution, then $X \sim \mu$ means the distribution of $X$ is equal to $\mu$. Recall also the notation $Y \cdot A$ for congruence transformations, given by (\ref{eq:congruence}). 
\begin{proposition}[Transformation properties]  \label{prop:transf}
  Let $Y$ be a random variable in $\mathcal{P}_{m\,}$. For all $A \in \mathrm{GL}(m)$,
 \begin{subequations}
 \begin{equation} \label{eq:translation}
 Y \sim G(\bar{Y},\sigma) \, \Longrightarrow \,  Y \cdot A \sim G(\bar{Y}\cdot A ,\sigma)
 \end{equation}
 Moreover, (recall $I \in \mathcal{P}_m$ is the $m \times m$ identity matrix),
 \begin{equation} \label{eq:inversion}
 Y \sim G(I,\sigma) \, \Longrightarrow \,  Y^{-1} \sim G(I ,\sigma)
 \end{equation}
 \end{subequations}
\end{proposition}
\textbf{Proof\,: } This follows from Propositions \ref{prop:metric} and \ref{prop:integral}, in \ref{subsec:geo}. To show (\ref{eq:translation}), let $\varphi : \mathcal{P}_m \rightarrow \mathbb{R}$ be a test function. If  $Y \sim G(\bar{Y},\sigma)$ and $Z = Y \cdot A$, then the expectation of $\varphi(Z)$ is given by
\begin{eqnarray}
\nonumber \int_{\mathcal{P}_m} \varphi (Y\cdot A) \, p(Y|\,\bar{Y},\sigma)\, dv(Y) \,= \, \int_{\mathcal{P}_m} \varphi (Z) \, p(Z\cdot A^{-1}|\,\bar{Y},\sigma)\, dv(Z)  
\end{eqnarray}
where the equality is a result of (\ref{eq:integration1}), and the variable of integration was simply renamed from $Y$ to $Z$. Now, by (\ref{eq:invmetric2}), $p(Z\cdot A^{-1}|\,\bar{Y},\sigma) = p(Z|\,\bar{Y}\cdot A,\sigma)$. This completes the proof of (\ref{eq:translation}). The proof of (\ref{eq:inversion}) follows a similar reasoning. \hfill $\blacksquare$ 
\vspace{0.1cm}
\begin{proposition}[Gaussian distribution via polar coordinates]  \label{prop:gaussianpolar}
Let $U$ and $r = (r_1, \ldots, r_m)$ be independent random variables, with values in $\mathrm{O}(m)$ and $\mathbb{R}^m$, respectively. Assume that $U$ is uniformly distributed on $\mathrm{O}(m)$ and that $r_1, \ldots, r_m$ have the following joint probability density
\begin{equation} \label{eq:eigendistribution}
p(r) =  (m!\, 2^m)^{\scriptscriptstyle -1} \, \times \omega_m \, \times \, 8^{\frac{m(m-1)}{4}} \times \,\zeta^{\scriptscriptstyle -1}(\sigma) \, e^{- (r^{\,2}_1 + \ldots + r^{\,2}_m)/2\sigma^2} \, \prod_{i < j } \sinh\left(|r_i - r_j|/2\right)
\end{equation}
If $Y = Y(r,U)$ as in (\ref{eq:polar}), then $Y \sim G(I,\sigma)$.
\end{proposition}
\vspace{0.2cm}
\textbf{Proof\,: } Assume $Y = Y(r,U)$ and let $\varphi : \mathcal{P}_m \rightarrow \mathbb{R}$ be a test function. It is enough to prove that the expectation of $\varphi(Y)$ is given by
\begin{equation} \label{eq:expecone}
  \int_{\mathcal{P}_m} \varphi(Y)\, p(Y|\, I,\sigma) \, dv(Y)
\end{equation}
Recall, (see~\cite{muirhead}, Page $70$), the uniform distribution on $\mathrm{O}(m)$ is given by $\omega^{\scriptscriptstyle -1}_m \times \det(\theta)$. Since $U$ and $r$ are independent, the expectation of $\varphi(Y)$ is given by
$$
\omega^{\scriptscriptstyle -1}_m \, \times \, \int_{\mathrm{O}(m)}\,\int_{\mathbb{R}^m}  \varphi(Y(r,U))\, \det(\theta) \, p(r) \prod^m_{i=1} dr_i
$$
Replacing the expression (\ref{eq:eigendistribution}) of $p(r)$, and rearranging, this is found equal to 
$$
(m!\, 2^m)^{\scriptscriptstyle -1} \,\times \, 8^{\frac{m(m-1)}{4}} \times \,\zeta^{\scriptscriptstyle -1}(\sigma) \, \int_{\mathrm{O}(m)}\,\int_{\mathbb{R}^m}  \varphi(Y(r,U))\, e^{- (r^{\,2}_1 + \ldots + r^{\,2}_m)/2\sigma^2} \, \det(\theta) \, \prod_{i < j } \sinh\left(|r_i - r_j|/2\right) \prod^m_{i=1} dr_i
$$
By (\ref{eq:integrationew}) and (\ref{eq:zproof1}), this is the same as  (\ref{eq:expecone}). The proof is thus complete. \hfill $\blacksquare$ \\[0.15cm]
\indent \textbf{Remark\,:} Proposition (\ref{prop:gaussianpolar}) has the following implication. If $Y \sim G(I,\sigma)$ then  \textit{the determinant of $Y$ has log-normal distribution}. Precisely, if $t = \log\det(Y)$ then $t$ is normally distributed with mean $0$ and variance $m \,\times\, \sigma^2$. This can be seen from (\ref{eq:eigendistribution}), using the fact that $t = r_1 + \ldots + r_{m\,}$. While interesting in itself, it is not used in this paper. \\[0.1cm]
\indent Propositions \ref{prop:transf} and \ref{prop:gaussianpolar} yield the following general method for sampling from a Gaussian distribution. \\[0.1cm]
--- \textbf{\textit{Sampling from $G(\bar{Y},\sigma)$}}\,: \\[0.1cm]
\phantom{s} $\blacktriangleright\,$ By Proposition \ref{prop:transf}, if $Y \sim G(I,\sigma)$ then $Y \cdot \bar{Y}^{1/2} \sim G(\bar{Y},\sigma)$. Therefore, it is enough to consider sampling from $G(I,\sigma)$. \\[0.1cm]
\phantom{s} $\blacktriangleright\,$ By Proposition \ref{prop:gaussianpolar}, to sample from $G(I,\sigma)$, it is enough to know how to sample from a uniform distribution on $\mathrm{O}(m)$, in order to generate $U$, and from the multivariate density (\ref{eq:eigendistribution}), in order to generate $r$. Once these are obtained, they can be replaced in the spectral decomposition $Y = Y(r,U)$, to generate $Y\!\sim\!G(I,\sigma)$. \\[0.1cm]
\phantom{s} $\blacktriangleright\,$ Sampling from a uniform distribution on $\mathrm{O}(m)$ is fairly straightforward, (see~\cite{muirhead}, Page $70$). Indeed, if $A$ is an $m \times m$ matrix, whose entries are independent and each with normal distribution $\mathcal{N}(0,1)$, and if $A = UT$ with $U$ orthogonal and $T$ upper triangular, then $U$ is uniformly distributed on $\mathrm{O}(m)$. \\[0.1cm]
\phantom{s} $\blacktriangleright\,$ Sampling from the density (\ref{eq:eigendistribution}) can be achieved by usual application of the Metropolis-Hastings algorithm~\cite{casella}. In the special case $m = 2$, the density (\ref{eq:eigendistribution}) can be written as a product of two densities, for two independent variables. Indeed, in this case, let $t = r_{\scriptscriptstyle 1} + r_{\scriptscriptstyle 2}$ and $\rho = r_{\scriptscriptstyle 1} - r_{\scriptscriptstyle 2\,}$. It follows from (\ref{eq:eigendistribution}) that 
\begin{equation} \label{eigendistributionm2}
p(r) \propto e^{- t^{\,2}/4\sigma^2}\times p\,(\rho) \hspace{1cm} p\, (\rho) \propto e^{- \rho^{\,2}/4\sigma^2} \sinh(|\rho|/2)
\end{equation}
Therefore, when $m = 2$, sampling from the probability density (\ref{eq:eigendistribution}) only requires sampling from the univariate probability densities of $t$ and $\rho\,$.       

\subsection{Statistical inference problems} \label{subsec:gaussianstat}
It is clear that Rao's distance is the main ingredient in the definition of Gaussian distributions. Accordingly, it is to be expected that Rao's distance play an essential role in statistical inference problems for Gaussian distributions. Here, this is considered for the problems of maximum likelihood estimation and of hypothesis testing. 
  
The following Proposition \ref{prop:mle} deals with maximum likelihood estimation. For a Gaussian distribution $G(\bar{Y},\sigma)$, this proposition
shows how maximum likelihood estimates of the parameters $\bar{Y}$ and $\sigma$ can be computed.

The more interesting case is that of $\bar{Y}$. It turns out that, based on independent samples $Y_1, \ldots, Y_{\scriptscriptstyle N}$ from $G(\bar{Y},\sigma)$, the maximum likelihood estimate of $\bar{Y}$ is equal to the empirical Riemannian centre of mass of $Y_1, \ldots, Y_{\scriptscriptstyle N\,}$. The concept of empirical Riemannian centre of mass is among the most widely used concepts in the applied literature~\cite{afsari}\cite{moakher1}. Proposition \ref{prop:mle} shows that Gaussian distributions provide a statistical foundation for this concept.  
     
The empirical Riemannian centre of mass of $Y_1, \ldots, Y_{\scriptscriptstyle N}$ is the unique global minimiser $\hat{Y}_{\scriptscriptstyle N}$ of $\mathcal{E}_{\scriptscriptstyle N} : \mathcal{P}_m \rightarrow \mathbb{R}_+$, 
\begin{equation} \label{eq:empirical}
 \mathcal{E}_{\scriptscriptstyle N}(Y) = \frac{1}{\strut N}\sum^{\scriptscriptstyle N}_{n=1} d^{\,2}(Y,Y_n)
\end{equation}
Comparing this to (\ref{eq:variance}) of \ref{subsec:geo}, it can be seen $\hat{Y}_{\scriptscriptstyle N}$ is the Riemannian centre of mass of the empirical distribution
$\pi_{\scriptscriptstyle N} = \frac{1}{\mathstrut N} (\delta_{Y_1} + \ldots + \delta_{Y_{\scriptscriptstyle N}})$, where $\delta_Y$ denotes the Dirac measure concentrated at $Y \in \mathcal{P}_{m\,}$. Thus, existence and uniqueness of $\hat{Y}_{\scriptscriptstyle N}$ follow from Proposition \ref{prop:curvature}.

\begin{proposition}[MLE and centre of mass]  \label{prop:mle}
  Let $Y_1, \ldots, Y_{\scriptscriptstyle N}$ be independent samples from a Gaussian distribution $G(\bar{Y},\sigma)$. Based on these samples, the maximum likelihood estimate
of the parameter $\bar{Y}$ is the empirical Riemannian centre of mass $\hat{Y}_{\scriptscriptstyle N}$ of $Y_1, \ldots, Y_{\scriptscriptstyle N}$. Moreover, the maximum likelihood estimate of the parameter $\sigma$ is the solution $\hat{\sigma}_{\scriptscriptstyle N}$ of the equation, 
\begin{equation} \label{eq:mlesigma}
   \sigma^3 \times \frac{d  }{\mathstrut d\sigma}\log \zeta(\sigma)= \mathcal{E}_{\scriptscriptstyle N}(\,\hat{Y}_{\scriptscriptstyle N})
\end{equation}
where the unknown is $\sigma$. Both $\hat{Y}_{\scriptscriptstyle N}$ and $\hat{\sigma}_{\scriptscriptstyle N}$ exist and are unique for any realisation of the samples $Y_1, \ldots, Y_{\scriptscriptstyle N\,}$.
\end{proposition}
\vspace{0.2cm}
\textbf{Proof\,: } Using (\ref{eq:gaussianpdf}), and the fact that $Y_1, \ldots, Y_{\scriptscriptstyle N}$ are independent, the log-likelihood function, of the parameters $\bar{Y}$ and $\sigma$, can be written
$$
\sum^{\scriptscriptstyle N}_{n=1} \log \, p(Y_n|\, \bar{Y},\sigma) = -N \log \zeta(\sigma) - \frac{1}{2\sigma^2} \sum^{\scriptscriptstyle N}_{n=1} d^{\,2}(\bar{Y},Y_n)  
$$
Since the first term on the right-hand side does not contain $\bar{Y}$, the maximum likelihood estimate of $\bar{Y}$ can be found by maximising the second term. This is equivalent to minimising the sum of squared distances, which also appears in (\ref{eq:empirical}). The corresponding global minimum is precisely $\hat{Y}_{\scriptscriptstyle N}$, whose existence and uniqueness follow from Proposition \ref{prop:curvature}, (compare to~\cite{vemuri}, Corollary 2.2., Page $598$). To find the maximum likelihood estimate of $\sigma$, introduce the new parameter $\eta = -1/\sigma^2$ and let $\psi(\eta) = \log \zeta(\sigma)$. Then, $\psi(\eta)$ is a strictly convex function, by application of H\"older's inequality to (\ref{eq:z}). Therefore, the log-likelihood  function, in the above expression, is a strictly concave function of $\eta$. Its maximum is found by solving the equation $\psi^\prime(\eta) = \mathcal{E}_{\scriptscriptstyle N}(\hat{Y}_{\scriptscriptstyle N})$, which is equivalent to (\ref{eq:mlesigma}), (here, the prime denotes derivation). Existence and uniqueness of the solution $\hat{\sigma}_{\scriptscriptstyle N}$  of (\ref{eq:mlesigma}) now follow. Indeed, $\psi^\prime(\eta)$ is a strictly increasing function, since $\psi(\eta)$ is strictly convex. Also, $\psi^\prime(\eta)$ tends to $0$ as $\eta \rightarrow -\infty$ and to $+\infty$ as $\eta \rightarrow 0$. This implies the equation $\psi^\prime(\eta) = c$ has a unique solution for any $c > 0$. \hfill $\blacksquare$ \\[0.1cm]
\indent Proposition \ref{prop:mle} indicates how to compute numerically the maximum likelihood estimates $\hat{Y}_{\scriptscriptstyle N}$ and $\hat{\sigma}_{\scriptscriptstyle N\,}$. Indeed, this proposition states that $\hat{Y}_{\scriptscriptstyle N}$ is the empirical Riemannian centre of mass of the samples $Y_1, \ldots, Y_{\scriptscriptstyle N\,}$. As such, the task of computing $\hat{Y}_{\scriptscriptstyle N}$ is well-studied in recent literature, and can be carried out using algorithms based on deterministic line-search~\cite{lenglet}\cite{ferreira}\,, or on stochastic gradient descent~\cite{bonnabel}\cite{arnaudonstoch}. The deterministic Riemannian gradient descent algorithm for computing $\hat{Y}_{\scriptscriptstyle N}$ is given below in formulae (\ref{eq:gradientdescent1}) and (\ref{eq:gradientdescent2}) of Paragraph \ref{subsec:em}, based on~\cite{lenglet}. With regard to $\hat{\sigma}_{\scriptscriptstyle N\,}$, Proposition \ref{prop:mle} states that it is the unique solution of equation (\ref{eq:mlesigma}). This is a non-linear equation in one real-valued unknown, namely $\sigma$. As such, its solution $\hat{\sigma}_{\scriptscriptstyle N}$ can be computed by standard application of Newton's algorithm.
 
Proposition \ref{prop:mle} also shows that $\hat{\sigma}_{\scriptscriptstyle N}$ measures the mean dispersion of the samples $Y_1, \ldots, Y_{\scriptscriptstyle N\,}$ away from $\hat{Y}_{N\,}$. Indeed, the left-hand side of equation (\ref{eq:mlesigma}) depends only on the unknown $\sigma$, while its right-hand side has the fixed value ${\scriptstyle N^{-1} \,\times\,\sum^N_{n=1}} \;  d^{\,2}(\hat{Y}_{\scriptscriptstyle N},Y_n)$, given by (\ref{eq:empirical}). Since Proposition \ref{prop:mle} states the solution $\hat{\sigma}_{\scriptscriptstyle N}$ of this equation exists and is unique for any value of its right-hand side, it follows that there exists a function $\Phi : \mathbb{R}_+ \rightarrow \mathbb{R}_+\,$, such that
\begin{equation} \label{eq:F}
  \hat{\sigma}_{\scriptscriptstyle N} = \, \Phi  (\, {\scriptstyle N^{-1} \,\times\,\sum^N_{n=1}} \;  d^{\,2}(\hat{Y}_{\scriptscriptstyle N},Y_n) )
\end{equation}
In view of equation (\ref{eq:mlesigma}), this means that $\Phi$ is the inverse function of $\, \sigma \mapsto \sigma^3 \times d \log \zeta(\sigma)/d\sigma$. The argument of the function $\Phi$ in (\ref{eq:F}) gives the mean dispersion of the samples $Y_1, \ldots, Y_{\scriptscriptstyle N\,}$ away from $\hat{Y}_{N\,}$. On the other hand, from the proof of Proposition \ref{prop:mle}, it can be seen that $\Phi$ is a strictly increasing function. Indeed, up to the change of parameter $\eta \mapsto \sigma$, this function $\Phi$ is the same as the inverse function of $\psi^\prime(\eta)$, which is strictly increasing. Thus, (\ref{eq:F}) shows that $\hat{\sigma}_{\scriptscriptstyle N}$ is a direct measure of the mean dispersion ${\scriptstyle N^{-1} \,\times\,\sum^N_{n=1}} \;  d^{\,2}(\hat{Y}_{\scriptscriptstyle N},Y_n)\,$, obtained by application of the strictly increasing function $\Phi$.

The following Proposition \ref{prop:hypothesis} deals with the problem of hypothesis testing. Precisely, it deals with the problem of testing against a given value $Y_{\scriptscriptstyle 0}$ of $\bar{Y}$, while  $\sigma$ remains unknown. This problem is described by the null and alternative hypotheses, (see~\cite{inference}, for general background), 
\begin{equation} \label{eq:hypotheses}
  H_{\scriptscriptstyle 0} : \bar{Y} = Y_{\scriptscriptstyle 0} \hspace{0.9cm}   H_{\scriptscriptstyle 1} : \bar{Y} \neq Y_{\scriptscriptstyle 0} 
\end{equation}
Proposition \ref{prop:hypothesis} expresses the log-likelihood ratio statistic $T$, used for testing $H_{\scriptscriptstyle 1}$ against $H_{\scriptscriptstyle 0\,}$.  
\vspace{0.1cm}
\begin{proposition}[Likelihood ratio test] \label{prop:hypothesis}
  The log-likelihood ratio statistic, corresponding to hypotheses (\ref{eq:hypotheses}), is the following
\begin{eqnarray} \label{eq:lrt}
T = \log \left( \frac{\zeta(\, \Phi(\,{\scriptstyle N^{-1} \,\times\, \sum^N_{n=1}} \;  d^{\,2}(\,\hat{Y}_{\scriptscriptstyle N}\,,\,Y_n) ))}{\strut \zeta(\, \Phi( \,{\scriptstyle N^{-1} \,\times\, \sum^N_{n=1}} \;  d^{\,2}(\,Y_{\scriptscriptstyle 0}\,,\,Y_n) ))} \right)  + 
 \frac{{\scriptstyle \sum^N_{n=1}} \;  d^{\,2}(\,\hat{Y}_{\scriptscriptstyle N}\,,\,Y_n) }{\strut 2 \Phi(\,{\scriptstyle N^{-1} \,\times\, \sum^N_{n=1}} \;  d^{\,2}(\,\hat{Y}_{\scriptscriptstyle N}\,,\,Y_n)) } \, - \, \frac{{\scriptstyle \sum^N_{n=1}} \;  d^{\,2}(\,Y_{\scriptscriptstyle 0}\,,\,Y_n) }{\strut 2 \Phi(\,{\scriptstyle N^{-1} \,\times\, \sum^N_{n=1}} \;  d^{\,2}(\,Y_{\scriptscriptstyle 0}\,,\,Y_n) )}
\end{eqnarray}
The likelihood ratio test rejects $H_{\scriptscriptstyle 0}$ when $T > k$ and accepts it otherwise, where $k$ is chosen to fix the level of significance.
\end{proposition}
\vspace{0.2cm}
\textbf{Proof\,: } By definition~\cite{inference}, the log-likelihood ratio statistic is
$$
T = \log \, \left\lbrace \, \sup \; {\scriptstyle \prod^{\scriptscriptstyle N}_{n=1}}\,  p(Y_n|\, \bar{Y},\sigma)\, \middle/ \, \sup_{H_{\scriptscriptstyle 0}} \; {\scriptstyle \prod^{\scriptscriptstyle N}_{n=1}}\,  p(Y_n|\, \bar{Y},\sigma) \, \right\rbrace
$$
where the supremum in the numerator is over all possible values of $(\bar{Y},\sigma)$, and the supremum in the denominator is over only those values which belong to hypothesis $H_{{\scriptscriptstyle 0}\,}$. Expression (\ref{eq:lrt}) follows by noting that the supremum in the numerator is realised for $(\bar{Y},\sigma) = (\hat{Y}_{\scriptscriptstyle N},\hat{\sigma}_{\scriptscriptstyle N})$, by definition of maximum likelihood estimates, and the supremum in the denominator is realised for $(\bar{Y},\sigma) = (Y_{\scriptscriptstyle 0},\hat{\sigma}_{\scriptscriptstyle 0})$, where 
$$
\hat{\sigma}_{\scriptscriptstyle 0} = \Phi(\,{\scriptstyle N^{-1} \,\times\, \sum^N_{n=1}} \;  d^{\,2}(\,Y_{\scriptscriptstyle 0}\,,\,Y_n) )
$$
The remaining part of the proposition is the usual definition of the likelihood ratio test, as in~\cite{inference}. \hfill $\blacksquare$ \\[0.1cm]
\indent Admittedly, expression (\ref{eq:lrt}) for the log-likelihood ratio $T$ is quite involved. This is due to the fact that hypotheses $H_{\scriptscriptstyle 0}$ and $H_{\scriptscriptstyle 1}$ of (\ref{eq:hypotheses}) are composite hypotheses, depending on the unknown parameter $\sigma$. If it were known that $\sigma = \sigma_{\scriptscriptstyle 0\,}$, then $T$ would reduce to, (as in the proof of Proposition \ref{prop:hypothesis}, this is by definition of the log-likelihood ratio),
\begin{equation} \label{eq:lrtsimple} 
    T \,= \, \frac{1}{\mathstrut 2\sigma_{\scriptscriptstyle 0}} \, \left \lbrace \, {\scriptstyle \sum^N_{n=1}} \;  d^{\,2}(\,\hat{Y}_{\scriptscriptstyle N}\,,\,Y_n) \, - \, {\scriptstyle \sum^N_{n=1}} \;  d^{\,2}(\,Y_{\scriptscriptstyle 0}\,,\,Y_n) \right \rbrace
\end{equation}
so the likelihood ratio test amounts to comparing the dispersion of the samples $Y_1,\ldots, Y_{\scriptscriptstyle N}$ away from $\hat{Y}_{\scriptscriptstyle N\,}$, to their dispersion away from $Y_{\scriptscriptstyle 0\,}$. In all generality, it can be shown on the basis of Proposition \ref{prop:normality} of Paragraph \ref{subsec:gaussianasymp}, that the asymptotic distribution of $T$ under hypothesis $H_{\scriptscriptstyle 0}$ is a chi-squared distribution (this is in the limit $N \rightarrow \infty$). Using this asymptotic distribution, it becomes straightforward to choose the value of $k$ in order to obtain a preassigned level of significance (\textit{i.e.} false alarm probability, or alpha-level). These statements yield a direct generalisation of Wilks' theorem~\cite{inference} (Page $132$).
\subsection{Riemannian centre of mass and asymptotic properties} \label{subsec:gaussianasymp}
The present paragraph recovers the asymptotic properties of the maximum likelihood estimate $\hat{Y}_{\scriptscriptstyle N}$ of the parameter $\bar{Y}$, in the limit $N \rightarrow \infty\,$. Recall, from Proposition \ref{prop:mle}, the estimate $\hat{Y}_{\scriptscriptstyle N}$ is the empirical Riemannian centre of mass of samples $Y_1, \ldots, Y_{\scriptscriptstyle N}$ from the distribution $G(\bar{Y},\sigma)$.

These asymptotic properties of consistency and of asymptotic normality of $\hat{Y}_{\scriptscriptstyle N}$ will be given in Propositions \ref{prop:consistency} and \ref{prop:normality}, respectively. Proposition \ref{prop:consistency} states that $\hat{Y}_{\scriptscriptstyle N}$ converges almost surely to $\bar{Y}$, while Proposition \ref{prop:normality} states that $\hat{Y}_{\scriptscriptstyle N\,}$, in a sense made precise by the proposition, is asymptotically normally distributed about $\bar{Y}$. 

It is outside the scope of the present paper to establish the further property of asymptotic efficiency of $\hat{Y}_{\scriptscriptstyle N\,}$. This property can be given a precise formulation on the basis of the work of Smith~\cite{smith05}, (see Corollary $3$, Page $1619$), and is the focus of ongoing work.

Propositions \ref{prop:consistency} and \ref{prop:normality} rely on the following Proposition \ref{prop:ybar}, which states that $\bar{Y}$ is the Riemannian centre of mass of the distribution $G(\bar{Y},\sigma)$, (in the sense defined after (\ref{eq:variance}) in \ref{subsec:geo}). 

In the statement of Proposition \ref{prop:ybar}, the following notation is used. Recall definition (\ref{eq:variance}) of the variance function $\mathcal{E}_\pi$ of a probability distribution $\pi$ on $\mathcal{P}_{m\,}$. If $\pi$ is the distribution $G(\bar{Y},\sigma)$, then $\mathcal{E}_\pi(Y)$ will be denoted $\mathcal{E}(Y| \, \bar{Y},\sigma)\,$.  
\vspace{0.1cm}
\begin{proposition}[Centre of mass of a Gaussian distribution]  \label{prop:ybar}
  For any $\bar{Y} \in \mathcal{P}_m$ and $\sigma > 0\,$, the following properties hold for the Gaussian distribution $G(\bar{Y},\sigma)$. \\[0.1cm]
(i) $\bar{Y}$ is the Riemannian centre of mass of $G(\bar{Y},\sigma)$. That is,
\begin{subequations}
\begin{equation} \label{eq:gaussbarycentre}
   \bar{Y} = \mathrm{argmin}_{ \scriptstyle Y} \; \mathcal{E}(Y| \, \bar{Y},\sigma)
\end{equation}
(ii) $\sigma$ is given by
\begin{equation} \label{eq:gaussigma}
   \sigma = \Phi  \left(\, {\textstyle \bigintsss_{\scriptscriptstyle \mathcal{P}_{m}}}\, d^{\,2}(\bar{Y},Z)\,p(Z|\, \bar{Y},\sigma) dv(Z)\, \right)
\end{equation}
where the function $\Phi$ was introduced in (\ref{eq:F}).
\end{subequations}
\end{proposition}
\vspace{0.2cm}
\textbf{Proof\,: }  The proof of this proposition is given in Appendix \ref{appendix0}. A more general, information-theoretic, proof of item (i) was given by Pennec~\cite{pennec1}, (Theorem 3, Page $142$). \hfill $\blacksquare$ \\[0.1cm]
\indent Proposition \ref{prop:consistency} states that the maximum likelihood estimate $\hat{Y}_{\scriptscriptstyle N}$ converges almost surely to the true value of the parameter $\bar{Y}$, when $N \rightarrow \infty\,$. 

This property of $\hat{Y}_{\scriptscriptstyle N}$ is called consistency. In asymptotic statistics, consistency of maximum likelihood estimates is a well-known general result~\cite{inference}. However, any statement of this result requires several additional technical conditions. Here, no such conditions will be necessary. Proposition \ref{prop:consistency} follows immediately from Proposition \ref{prop:ybar}, using a theorem on the convergence of empirical Riemannian centres of mass~\cite{bhatta1}. 
\vspace{0.1cm}
\begin{proposition}[Consistency of $\hat{Y}_{\scriptscriptstyle N\,}$] \label{prop:consistency}
 Let $Y_{\scriptscriptstyle 1}, Y_{\scriptscriptstyle 2}, \ldots$ be independent samples from a Gaussian distribution $G(\bar{Y},\sigma)$. The empirical Riemannian centre of mass $\hat{Y}_{\scriptscriptstyle N}$ of $Y_{\scriptscriptstyle 1},\ldots, Y_{\scriptscriptstyle N}$ converges almost surely to $\bar{Y}$, as $N \rightarrow \infty\,$. 
\end{proposition}
\vspace{0.1cm}
\textbf{Proof\,: } According to~\cite{bhatta1} (Theorem $2.3.$, Page $8$), if $Y_{\scriptscriptstyle 1}, Y_{\scriptscriptstyle 2}, \ldots$ are independent samples from a probability distribution $\pi$ on $\mathcal{P}_{m\,}$, and $\hat{Y}_{\scriptscriptstyle N}$ is the empirical Riemannian centre of mass of  $Y_{\scriptscriptstyle 1},\ldots, Y_{\scriptscriptstyle N}\,$, then $\hat{Y}_{\scriptscriptstyle N}$ converges almost surely to the Riemannian centre of mass of $\pi$, as $N \rightarrow \infty\,$. 

Applying this theorem to $\pi = G(\bar{Y},\sigma)$, it follows that $\hat{Y}_{\scriptscriptstyle N}$ converges almost surely to the Riemannian centre of mass of $G(\bar{Y},\sigma)$. However, by Proposition \ref{prop:ybar}, this Riemannian centre of mass is exactly $\bar{Y}$. The proposition is thus proved. \hfill $\blacksquare$ \\[0.12cm]
\indent Proposition \ref{prop:normality} is concerned with the asymptotic distribution of, roughly speaking, the ``difference" between $\hat{Y}_{\scriptscriptstyle N}$ and $\bar{Y}$. When the space $\mathcal{P}_m$ is considered as a Riemannian manifold, equipped with the Rao-Fisher metric, the natural definition of ``difference" between $\hat{Y}_{\scriptscriptstyle N}$ and $\bar{Y}$ is using the Riemannian logarithm mapping, (already mentioned in the proof of Proposition \ref{prop:ybar}, formula (\ref{eq:gradlog})). This will be denoted
\begin{equation} \label{eq:delta}
  \Delta = \mathrm{Log}_{\scriptscriptstyle \bar{Y}} (\hat{Y}_{\scriptscriptstyle N}) 
\end{equation}
The Riemannian logarithm mapping is defined as follows. Recall that any $Y,Z \in \mathcal{P}_m$ are connected by a unique geodesic curve $\gamma(t)$, given by  (\ref{eq:geodesic1}). The Riemannian logarithm $\mathrm{Log}_{\scriptscriptstyle Y}(Z)$ is an $m \times m$ symmetric matrix, equal to the initial derivative $\frac{d\gamma}{dt}(0)$. From (\ref{eq:geodesic1}), 
\begin{equation} \label{eq:log}
   \mathrm{Log}_{\scriptscriptstyle Y} (Z) = Y^{1/2} \log( Y^{-1/2} Z Y^{-1/2} ) \,Y^{1/2}
\end{equation}
where $\log$ is the symmetric matrix logarithm. 

In order to study the asymptotic distribution of the matrix $\Delta$ of (\ref{eq:delta}), (that is, the asymptotic joint distribution of matrix entries), it is desirable to have a vector representation, identifying $\Delta$ with an array $(\Delta_a; a = 1, \ldots, p)$ where each $\Delta_a$ is real-valued. 

To obtain this, consider for each $Y \in \mathcal{P}_m$ the scalar product, defined for any symmetric $m \times m$ matrices $v$ and $w$,
\begin{equation} \label{eq:scalarproduct}
 \langle v,w \rangle_{\scriptscriptstyle Y} = \mathrm{tr}\, [ Y^{\scriptscriptstyle -1}v \, Y^{\scriptscriptstyle  -1}w] 
\end{equation}
This scalar product is equivalent to the Rao-Fisher metric (\ref{eq:metric}), (in Riemannian geometry, this is taken to be the definition of the metric~\cite{helgason}\cite{terras2}). 

Let $\Delta_a$ be given by
\begin{equation} \label{eq:vector}
 \Delta_a  \, = \, \langle \Delta, e_a \rangle_{\scriptscriptstyle \bar{Y}}
\hspace{0.5cm} a = 1, \ldots, p
\end{equation} 
where $(e_a; a = 1, \ldots, p)$ is an orthonormal basis for the scalar product (\ref{eq:scalarproduct}), taken at $Y = \bar{Y}$. 

Then, $(\Delta_a;a = 1, \ldots, p)$ are the components of $\Delta$ in the basis $(e_a; a = 1, \ldots, p)$, so the asymptotic distribution of $\Delta$ is completely determined by the asymptotic distribution of the vector representation $(\Delta_a; a = 1, \ldots, p)$. 

Proposition \ref{prop:normality} states that the asymptotic joint distribution of the scaled components $N^{\scriptscriptstyle 1/2} \times \Delta_a$ is a centred normal distribution with covariance matrix equal to $4\,\sigma^{\scriptscriptstyle 4} \times C^{\scriptscriptstyle -1\,}$, where $C$ is the $p \times p$ symmetric positive definite matrix
\begin{equation} \label{eq:C}
  C_{ab} =  4 \times {\textstyle \bigintsss_{\scriptscriptstyle \mathcal{P}_{m}}}  
\left[\Delta_a(Z)\times \Delta_b(Z)\right] \,  p(Z|\, \bar{Y},\sigma) dv(Z) \hspace{0.5cm} a,b = 1, \ldots, p
\end{equation}
with, as in (\ref{eq:delta}) and (\ref{eq:vector}),
$$
\Delta_a(Z) = \langle\mathrm{Log}_{\scriptscriptstyle \bar{Y}} (Z),e_a\rangle_{\scriptscriptstyle \bar{Y}} \hspace{1cm} \mbox{ for } 
Z \in \mathcal{P}_{\scriptscriptstyle m}
$$ 

In the statement of Proposition \ref{prop:normality}, $\mathcal{L} \lbrace \cdot \rbrace$ denotes the probability distribution of the quantity in curly brackets, and $\Rightarrow$ denotes convergence of probability distributions. 
\vspace{0.1cm}
\begin{proposition}[Asymptotic normality of $\hat{Y}_{\scriptscriptstyle N}$] \label{prop:normality}
Let $(\Delta_a ; a = 1, \ldots, p)$ be given by (\ref{eq:vector}). In the limit $N \rightarrow \infty$, 
\begin{equation} \label{eq:normalityprop}
   \mathcal{L} \lbrace N^{\scriptscriptstyle 1/2} \times (\Delta_1,\ldots, \Delta_p) \rbrace \, \Rightarrow \mathcal{N}(0,4\,\sigma^{\scriptscriptstyle 4} \times C^{\scriptscriptstyle -1})
\end{equation}
where $\mathcal{N}(0,4\,\sigma^{\scriptscriptstyle 4} \times C^{\scriptscriptstyle -1})$ is the centred normal distribution on $\mathbb{R}^p$ with covariance matrix $4\,\sigma^{\scriptscriptstyle 4} \times C^{\scriptscriptstyle -1}$, with $C$ given by (\ref{eq:C}).
\end{proposition}
\vspace{0.1cm}
\textbf{Proof\,: } The proof follows from two lemmas, stated below. Lemma \ref{lemma:normality1} establishes asymptotic normality of $N^{\scriptscriptstyle 1/2} \times (\Delta_1,\ldots, \Delta_p)\,$, and Lemma \ref{lemma:normality2} shows the asymptotic covariance matrix is equal to $4\,\sigma^{\scriptscriptstyle 4} \times C^{\scriptscriptstyle -1\,}$.
\vspace{0.1cm}
\begin{lemma}[Asymptotic normality] \label{lemma:normality1}
  In the limit $N \rightarrow \infty$, 
\begin{equation} \label{eq:normality1}
   \mathcal{L} \lbrace N^{\scriptscriptstyle 1/2} \times (\Delta_1,\ldots, \Delta_p) \rbrace \, \Rightarrow \mathcal{N}(0,\Lambda)
\end{equation} 
where $\Lambda =   H^{-1} \, C \, H^{-1}$ and $H$ is the $p \times p$ symmetric positive definite matrix,
\begin{equation} \label{eq:H}
  H_{ab} = \nabla^2 \mathcal{E}(\bar{Y})(e_a,e_b) \hspace{0.5cm} a,b = 1, \ldots, p
\end{equation}
Here, $\mathcal{E}$ is the variance function $Y \mapsto \mathcal{E}(Y|\, \bar{Y},\sigma)$, defined before Proposition \ref{prop:ybar}, and  
$\nabla^2 \mathcal{E}(\bar{Y})$ is the Riemannian Hessian of this function, evaluated at $\bar{Y}$.
\end{lemma}
\vspace{0.1cm}
\begin{lemma}[Asymptotic covariance] \label{lemma:normality2}
The asymptotic covariance matrix $\Lambda$, appearing in (\ref{eq:normality1}), verifies $\Lambda = 4\,\sigma^{\scriptscriptstyle 4} \times C^{\scriptscriptstyle -1\,}$.
\end{lemma}
\vspace{0.1cm}

The proofs of these two lemmas are given in Appendix \ref{appendix}. A result similar to Lemma \ref{lemma:normality1} can be found in~\cite{bhatta2},
(Remark $2.2.$, Page $1232$), while Lemma \ref{lemma:normality2} is motivated by a result in~\cite{smith05}, (Theorem $1$, Page $1614$).

If Lemmas \ref{lemma:normality1} and \ref{lemma:normality2} are admitted, the proposition follows immediately. Indeed, replacing $\Lambda = 4\,\sigma^{\scriptscriptstyle 4} \times C^{\scriptscriptstyle -1\,}$, as stated in Lemma \ref{lemma:normality2}, into (\ref{eq:normality1}), produces (\ref{eq:normalityprop}). \hfill $\blacksquare$ 
\section{Mixtures of Riemannian Gaussian distributions} \label{sec:classif}
Mixtures of parameterised probability distributions on a Euclidean space are widely studied in the statistical literature~\cite{mixtures}\cite{mixture}. 
The present section generalises these mixtures, from Euclidean space to the space $\mathcal{P}_{m\,}$, considered as a Riemannian manifold. It does so by introducing \textit{mixtures of Riemannian Gaussian distributions} on $\mathcal{P}_{m\,}$. 

For shortness, mixtures of Riemannian Gaussian distributions will be called mixtures of Gaussian distributions. A mixture of Gaussian distributions is a probability distribution on $\mathcal{P}_{m}$ with probability density function
\begin{equation} \label{eq:mixture}
  p(Y) = \sum^{\scriptscriptstyle M}_{\mu \,= 1} \varpi_\mu \times p(Y | \,\bar{Y}_\mu,\sigma_\mu)
\end{equation}
Precisely, this is a probability density with respect to the Riemannian volume element of $\mathcal{P}_{m\,}$, given by (\ref{eq:dv}). Here, $\varpi_1, \ldots, \varpi_{\scriptscriptstyle M}$ are non-zero positive weights which satisfy  $\varpi_1 + \ldots + \varpi_{\scriptscriptstyle M} = 1$, and each density $p(Y | \,\bar{Y}_\mu,\sigma_\mu)$ is given by (\ref{eq:gaussianpdf}), with parameters $\bar{Y}_\mu \in \mathcal{P}_{m\,}$ and $\sigma_\mu > 0$. The number $M$ will be called the number of mixture components. 

The remainder of this section is devoted to the estimation of mixtures of Gaussian distributions, and to their use in the classification of data in $\mathcal{P}_{m\,}$. The results which are developed are partially based on a recent work by the authors~\cite{icip}. They are aimed to generalise
analogous results, for mixtures of parameterised distributions on a Euclidean space.

Paragraph \ref{subsec:em} derives a new EM (expectation-maximisation) algorithm, for computing maximum likelihood estimates of mixture parameters
$\vartheta = \lbrace (\varpi_\mu,\bar{Y}_\mu,\sigma_\mu) ; \mu = 1, \ldots, M \rbrace$. 

Paragraph \ref{subsec:laplace} proposes a new Bayes classification rule, based on mixtures of Gaussian distributions, for the classification of data in $\mathcal{P}_{m\,}$. 

The present section is motivated by the idea that the class of mixtures of Gaussian distributions is expected to be sufficiently rich, in order to represent the statistical distribution of data in $\mathcal{P}_m$ which arise in real-world applications. Roughly, this is because one expects that any probability density on $\mathcal{P}_{m}$ can be approximated to any required precision by a mixture of Gaussian distributions, provided the number $M$ of mixture components is suitably large. As stated in the introduction, experimental verification of this idea is still ongoing~\cite{gretsi}.

\subsection{A new EM algorithm for mixture estimation} \label{subsec:em}
Let $Y_1, \ldots, Y_{\scriptscriptstyle N}$ be independent samples, drawn from a mixture of Gaussian distributions, given by (\ref{eq:mixture}). This paragraph considers the task of computing maximum likelihood estimates of mixture parameters $\vartheta = \lbrace (\varpi_\mu,\bar{Y}_\mu,\sigma_\mu) ; \mu = 1, \ldots, M \rbrace$, based on the samples $Y_1, \ldots, Y_{\scriptscriptstyle N\,}$.

This task will be realised using a new EM (expectation-maximisation) algorithm, which generalises, to the context of the Riemannian geometry of the space $\mathcal{P}_{m\,}$, the currently existing EM algorithms used in the estimation of mixtures of parameterised distributions on a Euclidean space~\cite{mixtures}.

It is assumed, throughout the following, that the number of mixture components, denoted $M$ in (\ref{eq:mixture}), is known and fixed. The problem of determining a suitable $M$, in view of the samples $Y_1, \ldots, Y_{\scriptscriptstyle N\,}$, is known as order selection. This problem is not considered here, as it falls outside the scope of the methods used in the present paper. In principle, it could be tackled by application of existing general methods, such as those based on information criteria~\cite{leroux}.

In the remainder of this paragraph, the new EM algorithm will be derived, based on the general principle set forth in the founding paper~\cite{em}. This algorithm iteratively updates an approximation $\hat{\vartheta} = \lbrace (\hat{\varpi}_\mu,\hat{Y}_\mu,\hat{\sigma}_\mu) \rbrace\,$ of the maximum likelihood estimates of mixture parameters $\vartheta$, by repeated application of so-called E (expectation) and M (maximisation) steps. In order to specify these two steps, the following setting is needed.

To begin, assume the mixture model (\ref{eq:mixture}) was generated based on latent variables $L_1, \ldots, L_{\scriptscriptstyle N\,}$, defined as follows, (compare to~\cite{mixtures}\,, Page $84$). Let $L_1, \ldots, L_{\scriptscriptstyle N}$ be independent identically distributed random variables, which take the values $1, \ldots, M$, with respective probabilities, (here, $\mathbb{P}$ denotes an underlying probability measure),
\begin{equation} \label{eq:latents}
  \mathbb{P}\,(L_n = \mu) = \varpi_\mu 
\end{equation}
In order to obtain samples $Y_1, \ldots, Y_{\scriptscriptstyle N\,}$, from the mixture distribution (\ref{eq:mixture}), assume the couples $(L_n,Y_n)$ are
independent, and
\begin{equation} \label{eq:complete}
  \mathbb{P}\,(Y_n = Y | \,L_n = \mu ) = p(Y|\,\bar{Y}_\mu,\sigma_\mu) 
\end{equation}
Intuitively, $L_n$ is the membership label of $Y_{n\,}$, so $L_n = \mu$ indicates that $Y_n$ belongs to component $\mu$ of the mixture.
 
The distribution of $Y_1, \ldots, Y_{\scriptscriptstyle N}$ can be recovered from (\ref{eq:latents}) and (\ref{eq:complete}). Indeed,
\begin{eqnarray}
\nonumber \mathbb{P}\,(Y_n = Y) =\sum^{\scriptscriptstyle M}_{\mu \,= 1} \mathbb{P}\,(L_n = \mu) \times \mathbb{P}\,(Y_n = Y |\, L_n = \mu ) \\[-0.01cm]
\nonumber = \sum^{\scriptscriptstyle M}_{\mu \,= 1} \varpi_\mu \times p(Y | \,\bar{Y}_\mu,\sigma_\mu) \hspace{2.3cm}
\end{eqnarray}
which is exactly (\ref{eq:mixture}).

Let $C_\mu$ denote the number of $L_n$ which are equal to $\mu\,$. That is,
\begin{equation} \label{eq:count}
  C_\mu = \sum^{\scriptscriptstyle N}_{n=1} \mathds{1}_{ \lbrace L_n = \, \mu \rbrace}
\end{equation}
where $\mathds{1}_{\scriptscriptstyle A}$ denotes the indicator function of event $A$, (equal to $1$ if $A$ holds and to zero otherwise). Given the interpretation of the $L_n$ as membership labels, $C_\mu$ is called the membership count of component $\mu$ of the mixture.

Consider now the quantities $\omega_\mu(Y_n)$ and $N_{\mu\,}$, defined as follows, (see discussion in~\cite{em}, Section 4.3.),
\begin{eqnarray} \label{eq:memberprob}
\nonumber \omega_\mu(Y_n) = \mathbb{E}\,( \mathds{1}_{ \lbrace L_n = \, \mu \rbrace} | \, Y_1, \ldots, Y_{\scriptscriptstyle N}) \\[0.1cm]
\phantom{\omega_\mu(Y} = \mathbb{P}\,( L_n = \, \mu | \, Y_1, \ldots, Y_{\scriptscriptstyle N}) \,\,\, 
\end{eqnarray}
\begin{equation} \label{eq:membercount}
N_\mu = \mathbb{E}\,[ C_\mu | \, Y_1, \ldots, Y_{\scriptscriptstyle N}] \hspace{2.2cm} \,
\end{equation}
where $\mathbb{E}$ denotes expectation. Given these definitions, one calls $\omega_\mu(Y_n)$ the conditional membership probability, and $N_\mu$ the conditional membership count, of component $\mu$ of the mixture. These are given by the following formulae,
\begin{equation} \label{eq:emdescription1}
   \omega_\mu(Y_n) \propto\varpi_\mu \times p(Y_n|\, \bar{Y}_\mu,\sigma_\mu) \hspace{0.6cm} N_\mu = \sum^{\scriptscriptstyle N}_{n=1} \omega_\mu(Y_n)
\end{equation}
which follows from (\ref{eq:memberprob}) and (\ref{eq:membercount}), after a simple application of Bayes' theorem. Here, the constant of proportionality, corresponding to ``$\propto$", is chosen so that $\omega_1(Y_n) + \ldots + \omega_{\scriptscriptstyle M}(Y_n)= 1$. Moreover, it is clear that
\begin{equation} \label{eq:emdescription2}
\sum^{\scriptscriptstyle M}_{\mu = 1} N_\mu = N
\end{equation}
The E and M steps can now be specified, using definitions (\ref{eq:memberprob}) and (\ref{eq:membercount}). To emphasise that these two definitions are applied under a given value of $\vartheta = \lbrace (\varpi_\mu,\bar{Y}_\mu,\sigma_\mu)\rbrace$, the following notation is used,
\begin{equation} \label{eq:emdescription3}
\omega_\mu(Y_n) = \omega_\mu(Y_n,\vartheta) \hspace{0.6cm} N_\mu = N_\mu(\vartheta)
\end{equation}
The general form of the E and M steps is the following~\cite{em}.\\[0.1cm]
\noindent $\blacktriangleright$ E step\,: based on the current value of $\hat{\vartheta}$, compute
\begin{equation} \label{eq:e}
  Q(\vartheta|\,\hat{\vartheta}\,) =  \, \mathbb{E}_{\,\hat{\vartheta}} \, \left[ \,{\scriptstyle \sum^N_{n=1}}\, \ell(L_n,Y_n) \, \middle|\, Y_1, \ldots, Y_{\scriptscriptstyle N\,}\right] 
\end{equation}
where $\mathbb{E}_{\,\hat{\vartheta}}$ denotes expectation under the value $\hat{\vartheta}$ of the mixture parameters $\vartheta$, and where $\ell(L_n,Y_n)$ is the joint log-likelihood of $L_n$ and  $Y_{n\,}$. In other words, the E step computes the conditional expectation of the complete log-likelihood, based on the current value of $\hat{\vartheta}$. \\[0.1cm]
\noindent $\blacktriangleright$ M step\,: assign to $\hat{\vartheta}$ the new value
\begin{equation} \label{eq:m}
 \hat{\vartheta}^{\mathrm{new}} = \mathrm{argmax}_{\,\vartheta} \;\, Q(\vartheta|\,\hat{\vartheta}\,) 
\end{equation}
In other words, the M step updates $\hat{\vartheta}$ by maximising the conditional expectation of the complete log-likelihood, as provided by the E step. \\[0.1cm]
\indent In order to carry out the E step, it is enough to note
\begin{equation}
 \ell(L_n,Y_n) = {\scriptstyle \sum^M_{\mu=1}} \, \mathds{1}_{ \lbrace L_n = \, \mu \rbrace}
\lbrace \log \, \varpi_\mu \,+\, \log \, p(Y_n|\, \bar{Y}_\mu,\sigma_\mu) \rbrace
\end{equation}
as follows from (\ref{eq:latents}) and (\ref{eq:complete}). After replacing expression (\ref{eq:gaussianpdf}) of $p(Y_n|\, \bar{Y}_\mu,\sigma_\mu)$ and performing some algebraic transformations, the conditional expectation in (\ref{eq:e}) is found to be, (see again~\cite{em}, Section 4.3.),
\begin{equation} \label{eq:q}
Q(\vartheta|\,\hat{\vartheta}\,) = {\scriptstyle \sum^M_{\mu=1}} \, \left \lbrace \, N_\mu(\hat{\vartheta})\, \lbrace \log \, \varpi_\mu - \log \, \zeta(\sigma_\mu) \rbrace \, - \,{\scriptstyle \sum^N_{n=1}} \,  \omega_\mu(Y_n,\hat{\vartheta}) \, \left. d^{\,2}(\bar{Y}_\mu,Y_n)\middle/ 2\sigma^2_\mu \right. \right \rbrace
\end{equation} 
\indent In order to carry out the M step, it is required to maximise this expression with respect to $\vartheta = \lbrace (\varpi_\mu,\bar{Y}_\mu,\sigma_\mu)\rbrace$. The maximisation can be carried out, first over the values of $\varpi_{\mu\,}$, then over those of $\bar{Y}_{\mu\,}$, and finally over those of $\sigma_{\mu\,}$.

Note that $\varpi_{\mu}$ enters (\ref{eq:q}) only via the expression
\begin{equation} \label{eq:qomega}
   \sum^M_{\mu=1} \, N_\mu(\hat{\vartheta})\,\log \, \varpi_\mu
\end{equation}
By Jensen's inequality, applied to the concave function $\log\,$, this is maximised, (under the constraint $\varpi_1 + \ldots + \varpi_{\scriptscriptstyle M} = 1$),
when $\varpi_\mu = \hat{\varpi}^{\mathrm{new}}_\mu$, where 
\begin{equation} \label{eq:omegaupdate}
  \hat{\varpi}^{\mathrm{new}}_{\mu} = \left. N_\mu(\hat{\vartheta}) \middle/ N  \right.
\end{equation}
Maximising (\ref{eq:q}) over the values of $\bar{Y}_{\mu\,}$ reduces to minimising, separately, each of the expressions,
\begin{equation} \label{eq:ybarupdate}
   \mathcal{E}_\mu(Y) = \sum^N_{n=1} \,  \omega_\mu(Y_n,\hat{\vartheta}) \, d^{\,2}(Y,Y_n)
\end{equation}
over the values of $Y \in \mathcal{P}_{m\,}$. The corresponding minima will be denoted $\hat{Y}^{\mathrm{new}}_{\mu}$.  

Finally, maximising (\ref{eq:q}) over the values of $\sigma_\mu$ can be carried out by differentiating with respect to $\sigma_{\mu\,}$, and setting the derivative equal to zero. A direct calculation shows that this yields the solution
\begin{equation} \label{eq:sigmaupdate}
   \hat{\sigma}^{\mathrm{new}}_{\mu} = \Phi  (\, {\scriptstyle N^{\scriptscriptstyle -1}_\mu(\hat{\vartheta}) \, \times \, \sum^N_{n=1}} \;  \omega_\mu(Y_n,\hat{\vartheta}) \, d^{\,2}(\hat{Y}_\mu,Y_n) )
\end{equation}
where the function $\Phi$ was defined in (\ref{eq:F}), Paragraph \ref{subsec:gaussianstat}.

From the above, it is seen that the new EM algorithm, from a practical point of view, consists in repeated application of the update rules (\ref{eq:omegaupdate}), (\ref{eq:ybarupdate}) and (\ref{eq:sigmaupdate}), in this same order. These update rules should be repeated for as long as they introduce a sensible change in the values of $\hat{\varpi}_{\mu\,}$, $\hat{Y}_{\mu\,}$, and $\hat{\sigma}_{\mu\,}$. Other stopping criteria, involving the amount of increase in the joint likelihood function of $\hat{\varpi}_{\mu\,}$, $\hat{Y}_{\mu\,}$, and $\hat{\sigma}_\mu$ can also be used.

Realisation of the update rules for $\hat{\varpi}_{\mu\,}$ and  $\hat{\sigma}_{\mu}$ is rather straightforward. On the other hand, the update rule for $\hat{Y}_\mu$ requires minimisation of the function $\mathcal{E}_\mu : \mathcal{P}_m \rightarrow \mathbb{R}_+$, defined by (\ref{eq:ybarupdate}). 
This is a function of the form given by (\ref{eq:variance}), in Paragraph \ref{subsec:geo}, so Proposition \ref{prop:curvature} guarantees that it has a unique global minimiser $\hat{Y}^{\mathrm{new}}_{\mu}$, which is also a unique stationary point.

Comparing (\ref{eq:ybarupdate}) to the general expression (\ref{eq:variance}), it appears clearly that $\hat{Y}^{\mathrm{new}}_{\mu}$ is the Riemannian centre of mass of the probability distribution $\omega_{\mu}$ on $\mathcal{P}_{m\,}$, given by
\begin{equation} \label{eq:omega}
  \omega_{\mu\,} = \sum^{\scriptscriptstyle N}_{n=1} \omega_\mu(Y_n) \times \delta_{Y_n}
\end{equation}
Here, recall that $\delta_Y$ denotes the Dirac measure concentrated at $Y \in \mathcal{P}_{m\,}$. As mentioned after Proposition \ref{prop:mle}, there exist several algorithms for computing Riemannian centres of mass~\cite{lenglet}\cite{ferreira,bonnabel, arnaudonstoch}. One of the most familiar among them is the Riemannian gradient descent algorithm, which may be found in~\cite{lenglet}. The $k$th iteration of this algorithm produces an approximation $\hat{Y}^{\scriptscriptstyle^k}_\mu$ of $\hat{Y}^{\mathrm{new}}_{\mu}$ in the following way. For $k = 1, 2, \ldots,$ let $\Delta_{\scriptscriptstyle k}$ be the symmetric matrix
\begin{equation} \label{eq:gradientdescent1}
  \Delta_{\scriptscriptstyle k} = \sum^{\scriptscriptstyle N}_{n=1} \omega_\mu(Y_n) \times \, \mathrm{Log}_{\,\hat{Y}^{\scriptscriptstyle^{k-1}}_\mu} \, (\,Y_n)
\end{equation}
where $\mathrm{Log}$ denotes the Riemannian logarithm mapping, given by (\ref{eq:log}). Then, $\hat{Y}^{\scriptscriptstyle^k}_\mu$ is defined to be
\begin{equation} \label{eq:gradientdescent2}
    \hat{Y}^{\scriptscriptstyle^k}_\mu = \mathrm{Exp}_{\,\hat{Y}^{\scriptscriptstyle^{k-1}}_\mu}\,  (\tau_{\scriptscriptstyle k} \, \Delta_{\scriptscriptstyle k})
\end{equation}
where  $\mathrm{Exp}$ is the Riemannian exponential mapping, inverse to the Riemannian logarithm mapping, given by
\begin{equation} \label{eq:exp}
  \mathrm{Exp}_Y \,(\Delta) = Y^{1/2}\, \exp\left(Y^{-1/2}\, \Delta \,Y^{-1/2}\right)\, Y^{1/2}
\end{equation}
with $\exp$ the matrix exponential, and where $\tau_{\scriptscriptstyle k} > 0$ is a step size, to be determined using a backtracking procedure.

The Riemannian gradient descent algorithm, specified by (\ref{eq:gradientdescent1}) and (\ref{eq:gradientdescent2}), is repeated as long as $\Vert \Delta_{\scriptscriptstyle k}\Vert > \epsilon\,$, where $\Vert \Delta_{\scriptscriptstyle k}\Vert$ is given by (\ref{eq:metric}) and $\epsilon$ is a precision
parameter, to be chosen by the user. This algorithm is guaranteed to converge, when a suitable backtracking procedure is used, regardless of the initialisation $\hat{Y}^{\scriptscriptstyle^0}_{\mu\,}$. A discussion of this convergence, in the case where Armijo backtracking procedure is used, can be found in~\cite{absil} (Theorem $4.3.1.$, Page $65$).

The Riemannian exponential mapping (\ref{eq:exp}), which is difficult to compute, can be replaced by certain so-called retraction mappings, which are less computationally demanding, without any change in the convergence of the algorithm. Some of these retraction mappings are given in~\cite{sra}, while a general characterisation of the class of suitable retraction mappings can be found in~\cite{absmal} (see, in particular, Paragraph $3.4.$, Page $11$).

\subsection{Classification using mixtures of Gaussian distributions} \label{subsec:laplace}
Mixture distributions lead to popular tools for classification of data which lie in a Euclidean space~\cite{tibshirani}. Moving beyond this usual Euclidean setting, the present paragraph considers the use of mixtures of Gaussian distributions in the classification of data in $\mathcal{P}_{m\,}$. 

Classification of data which lie in $\mathcal{P}_{m\,}$ is an important challenge for several applications, including remote
sensing~\cite{dong}\cite{wishartclass}, computer vision~\cite{tuzel}\cite{saintjean}, and medical imaging~\cite{congedo}\cite{medmix}. Most of these applications have used already existing classification techniques, such as Bayes classification (see~\cite{tibshirani}, Section $2.4.$), used in~\cite{wishartclass}\cite{saintjean}\cite{medmix}, or regression techniques (see~\cite{tibshirani}, Section $2.7.$), used in~\cite{dong}\cite{tuzel}. On the other hand, some applications have evolved new classification techniques based on the Riemannian geometry of $\mathcal{P}_{m\,}$. This is the case of~\cite{congedo}, whose approach, (discussed in the introduction of the present paper), is based on using Rao's Riemannian distance. 

In the present paragraph, the focus will be on supervised classification. A new classification rule is introduced which, like those of~\cite{wishartclass}\cite{saintjean}\cite{medmix}, implements the principle of Bayes classification, but which also uses Rao's Riemannian distance, as in~\cite{congedo}. Precisely, this new classification rule is a Bayes optimal classification rule, using posterior membership probabilities, computed by the EM algorithm, described in the previous paragraph.

The following setting is considered, (similar to~\cite{congedo}\cite{wishartclass}\cite{saintjean}\cite{medmix}). Assume known a training sequence $\mathcal{T}$. Precisely, $\mathcal{T} \subset \mathcal{P}_m$ is a set of data points, produced in some real-world application. Assume also known a partition of $\mathcal{T}$ into disjoint classes. Each class $\mathcal{C} \subset \mathcal{T}$ results from a well-defined, distinct experiment, (for example, application of a  measurement device to a given object). The data points which make up each class $\mathcal{C}$, while arising from the same experiment, may still display heterogeneous properties, (for example, if they are produced from measurements taken in different conditions). Therefore, each class $\mathcal{C}$ may further break down into clusters $C_1, \ldots, C_{\scriptscriptstyle M\,}$, where $M$ depends on the given class $\mathcal{C}$. Accordingly, the training sequence $\mathcal{T}$ can be partitioned directly into disjoint clusters, say $C_1, \ldots, C_{\scriptscriptstyle K\,}$, where $K$ is the total number of clusters within the training sequence. Then, the data points which belong to each cluster display essentially homogeneous properties.

In this setting, one has to carry out the two following tasks. First, for each class $\mathcal{C} \subset \mathcal{T}$, to identify the clusters $C_1, \ldots, C_{\scriptscriptstyle M\,}$ within this class. Second, whenever a new data point is produced, to associate this data point to a suitable cluster among $C_1, \ldots, C_{\scriptscriptstyle K\,}$. 

For the first task, it is required to perform a clustering analysis of the data points in each class $\mathcal{C}$. In~\cite{nielsen}\cite{saintjean}, these data points are modeled as a realisation of a mixture of Wishart distributions. Then, clustering analysis is performed using an EM algorithm for estimation of mixtures of Wishart distributions, or some variant of such an algorithm. 

Here, in view of introducing a classification rule which uses Rao's distance, Gaussian distributions are chosen over Wishart distributions. That is, the data points in each class $\mathcal{C}$ are modeled as a realisation of a mixture of Gaussian distributions, given by (\ref{eq:mixture}), rather than a mixture of Wishart distributions, as in~\cite{nielsen}\cite{saintjean}. Then, clustering analysis is performed using the EM algorithm of Paragrah \ref{subsec:em}.

Precisely, this algorithm computes $M$ triples of maximum likelihood estimates, $\hat{\vartheta} = \lbrace (\hat{\varpi}_\mu,\hat{Y}_\mu,\hat{\sigma}_\mu) ; \mu = 1, \ldots, M \rbrace\,$. These are used to identify clusters $\lbrace C_{\mu\,}; \mu = 1, \ldots, M\rbrace$ within the class $\mathcal{C}$, in the following way. If $\mathcal{C} = \lbrace Y_{n\,}; n = 1, \ldots, N \rbrace\,$, then each data point $Y_n$ is associated to the cluster $C_{\mu^*}$ which realises the maximum conditional membership probability $\omega_{\mu^*}(Y_n,\hat{\vartheta}) = \mathrm{max}_\mu \; \omega_\mu(Y_n,\hat{\vartheta})\,$, (recall the notation of (\ref{eq:emdescription1}) and (\ref{eq:emdescription3})).

With the first task being realised as just described, the training sequence $\mathcal{T}$ is partitioned into disjoint clusters, say $\lbrace C_{\kappa\,}; \kappa = 1,\ldots, K\rbrace$. Moreover, through the mixture model (\ref{eq:mixture}), each cluster $C_\kappa$ is identified with a Gaussian distribution $G(\bar{Y}(\kappa), \sigma(\kappa))\,$, and may be represented by a triple of maximum likelihood estimates, $(\hat{\varpi}(\kappa),\hat{Y}(\kappa),\hat{\sigma}(\kappa))$. 

For the second task, assume a new data point $Y_{t\,}$, called a test data point, has become available. The optimal way, (in the sense explained in~\cite{tibshirani}, Section $4.3.$), of associating this data point to a cluster $C_{\kappa\,}$, consists in applying the following Bayes classification rule\,: associate $Y_t$ to the cluster $C_{\kappa^*}$ which realises the maximum
\begin{equation} \label{eq:bayesclass}
N(\kappa^*) \times p(Y_t|\, C_{\kappa^*}) \, =\,   \mathrm{max}_{\kappa} \; N(\kappa) \times p(Y_t|\, C_\kappa) 
\end{equation} 
where $N(\kappa)$ is the number of data points in $C_{\kappa\,}$, and $p(Y_t|\,C_\kappa)$ denotes the density of $Y_{t\,}$, assuming it belongs to $C_{\kappa\,}$.  This is called a Bayes classification rule because the product $N(\kappa) \times p(Y_t|\, C_\kappa)$ implements Bayes formula, with $C_\kappa$ assigned a prior probability $\mathbb{P}(\kappa) \propto N(\kappa)$.

The new classification rule, introduced in the present paragraph, evaluates (\ref{eq:bayesclass}) using the representation of cluster $C_\kappa$ by  maximum likelihood estimates $(\hat{\varpi}(\kappa),\hat{Y}(\kappa),\hat{\sigma}(\kappa))$. Recall these are given by (\ref{eq:omegaupdate}), (\ref{eq:ybarupdate}) and (\ref{eq:sigmaupdate}). From (\ref{eq:omegaupdate}), $\hat{\varpi}(\kappa)$ is a consistent estimate of $\mathbb{P}(\kappa)$. On the other hand, in the present context, the density $p(Y_t|\,C_\kappa)$ is $p(Y_t| \,\bar{Y}(\kappa),\sigma(\kappa))$, which is consistently estimated using $p(Y_t| \,\hat{Y}(\kappa),\hat{\sigma}(\kappa))$. Accordingly, the new classification rule is the following\,: associate $Y_t$ to the cluster $C_{\kappa^*}$ which realises the maximum $\mathrm{max}_{\kappa} \, \hat{\varpi}(\kappa) \times p(Y_t|\, \hat{Y}(\kappa),\hat{\sigma}(\kappa))$.

Recalling expression (\ref{eq:gaussianpdf}) for $p(Y_t|\, \hat{Y}(\kappa),\hat{\sigma}(\kappa))$, the proposed classification rule takes on its final form. Precisely, $C_{\kappa^*}$ is the cluster which realises the minimum
\begin{equation} \label{eq:mlcrbis}
  \min_\kappa \, \left \lbrace - \log \, \hat{\varpi}(\kappa) + \log \, \zeta(\hat{\sigma}(\kappa)) + \frac{d^{\,2}(Y_{t\,},\hat{Y}(\kappa))}{\strut 2\hat{\sigma}^2(\kappa)}  \right \rbrace
\end{equation}
as can be found from (\ref{eq:gaussianpdf}) and (\ref{eq:bayesclass}) after taking logarithms and changing sign. 

The classification rule (\ref{eq:mlcrbis}) can be interpreted as follows. This rule prefers clusters $C_\kappa$ having a larger number of data points, (the minimum contains $-\log \, \hat{\varpi}(\kappa)$), or a smaller dispersion away from their Riemannian centre of mass, (the minimum contains $\log \, \zeta(\hat{\sigma}(\kappa))$). When choosing between two clusters with the same number of points and the same dispersion, this rule prefers the one whose Riemannian centre of mass is closer to $Y_{t\,}$. 

If the role of respective sizes and dispersions of clusters is neglected, (recall the size of a cluster refers to the number of its data points), then (\ref{eq:mlcrbis}) reduces to a nearest neighbour rule, which chooses $C_{\kappa^*}$ in order to realise the minimum
\begin{equation} \label{eq:neighborbis}
      \mathrm{min}_\kappa \, \lbrace d(Y_{t\,},\hat{Y}(\kappa)) \rbrace
\end{equation}
This nearest neighbour rule is the main subject of~\cite{congedo}.   

To close this paragraph, note that the general Bayes classification rule (\ref{eq:bayesclass}) was also applied in~\cite{wishartclass}, but with $p(Y_t|\, C_\kappa)$ given by the density of a Wishart distribution. Assume this Wishart distribution has expectation $\Sigma(\kappa) \in \mathcal{P}_m$ and number of degrees of freedom $n(\kappa)$, (for the definition of these parameters, see~\cite{muirhead}). Then, by the same reasoning leading from (\ref{eq:bayesclass}) to (\ref{eq:mlcrbis}), one obtains the Wishart classifier, which requires $C_{\kappa^*}$ to realise the minimum 
\begin{eqnarray} \label{eq:wishartclass}
\min_\kappa \, \left \lbrace - 2\log \, \hat{\varpi}(\kappa) - 
\hat{n}(\kappa) \, (\log\det \, (\hat{\Sigma}^{\scriptscriptstyle -1}(\kappa)Y_t) -  \mathrm{tr}(\hat{\Sigma}^{\scriptscriptstyle  -1}(\kappa)Y_t)  )  \right \rbrace
\end{eqnarray}
where $\hat{\varpi}(\kappa)$, $\hat{\Sigma}(\kappa)$ and $\hat{n}(\kappa)$ denote maximum likelihood estimates, which can be computed as in~\cite{nielsen}\cite{saintjean}.

The following Section compares the performance of classification rules (\ref{eq:mlcrbis}), (\ref{eq:neighborbis}) and (\ref{eq:wishartclass}), through a numerical experiment carried out on real data.

\section{Numerical experiment} \label{sec:numerical}
Among the many applications which require an effective approach to the classification of data in $\mathcal{P}_{m\,}$, is the problem of \textit{texture classification}, in computer vision. 

The present section applies the approach developed in Paragraph \ref{subsec:laplace}, to this problem. It shows that, in the context of this concrete application, the new classification rule (\ref{eq:mlcrbis}) offers significantly better performance than classification rules (\ref{eq:neighborbis}) and (\ref{eq:wishartclass}), proposed in \cite{congedo} and \cite{wishartclass}, respectively.

This section is organised as follows. First, a brief introduction to the problem of texture classification is given. Second, a numerical experiment is described, where classification rules (\ref{eq:mlcrbis}), (\ref{eq:neighborbis}) and (\ref{eq:wishartclass}) are applied to this problem. Finally, the results of this numerical experiment are summarised and commented, based on Table \ref{table:experiment}, below.

In computer vision~\cite{maria}, a texture is a pattern of local variation in image intensity, observed on a fixed scale. For example, text printed on white paper, as on the current page, constitutes a texture. In a remote sensing image, regions belonging to agricultural land, forest, and urban area give rise to different textures. 

The problem of texture classification involves textures belonging to one of a set of predefined classes. In remote sensing, these could be agricultural land, forest, urban area, \textit{etc}. The problem is to assign any new image, (more often, subregion of an image), to a suitable class. This problem is fundamental to medical image analysis, remote sensing, and material science, among other fields. 

The relationship between the problem of texture classification, on the one hand, and the problem of classification of data in $\mathcal{P}_{m\,}$, on the other hand, can be described as follows.

Many popular mathematical representations of texture are based on statistical modeling of wavelet coefficients~\cite{do,choy,verd,lasmar}. These representations have often been found useful, and are also justified by physiological and psychological studies of human visual
perception~\cite{psycho1}\cite{psycho2}. 

In recent work~\cite{verd}\cite{lasmar}, the statistical modeling of wavelet coefficients was carried out using multivariate probability distributions, which depend on a parameter $Y \in \mathcal{P}_{m\,}$. Such distributions include, for example, multivariate generalised Gaussian distributions, used in~\cite{verd}, which are parameterised by two positive scalar parameters, called the scale and shape parameters, and by $Y \in \mathcal{P}_{m\,}$, called the scatter matrix.

Accepting the mathematical representation that ``\textit{texture $\approx$ probability distribution of wavelet coefficients}," it is clear that, based on models such as those of~\cite{verd}\cite{lasmar}, a texture can be described by a parameter $Y \in \mathcal{P}_{m\,}$, (of course, one may wish to include other parameters). Then, the problem of texture classification becomes identical to the problem of classification of data in $\mathcal{P}_{m\,}$, described in \ref{subsec:laplace}.

In order to apply the classification rules considered in \ref{subsec:laplace} to the problem of texture classification, a numerical experiment was carried out, using the Vision Texture image database (VisTex)~\cite{vistex}. The Vision Texture database contains images, (size $512 \times 512$ pixels), of different materials or objects, such as bark, metal, brick, buildings, clouds, \textit{etc}. 

Starting with $40$ images, from this database, each image was subdivided into patches, (size $128 \times 128$ pixels, with $32$ pixel overlap). These patches are known to have the appropriate scale, at which textures can be observed. Therefore, it is possible to consider that each patch constitutes an individual texture, (This setup is usual in experiments using the VisTex database~\cite{do}\cite{lasmar}).

The subdivision of $40$ images into patches produces $169$ patches per image. Out of these patches, $84$ were used for training, and the remaining  $85$ for classification. Training refers to computation of the maximum likelihood estimates appearing in the classification rules (\ref{eq:mlcrbis}), (\ref{eq:neighborbis}) and (\ref{eq:wishartclass}). Once this is realised, classification is carried out by applying the rules (\ref{eq:mlcrbis}), (\ref{eq:neighborbis}) and (\ref{eq:wishartclass}) to the ``classification patches,'' in the aim of assigning each patch to a suitable class, in terms of the material or object which it represents.

Table \ref{table:experiment} shows the performance of each one of these three classification rules, in terms of so-called overall accuracy. Overall accuracy is the percentage of classification patches, (out of a total number of $85 \times 40$), which are correctly classified. The results given in Table \ref{table:experiment} are averaged over $100$ realisations of the experiment. For each realisation, the choice of training and classification patches, within the $169$ patches in each image, is generated at random. The mean and standard deviation of overall accuracy for each classification rule are given in the format $\mbox{\textit{mean} } \pm \mbox{ \textit{standard deviation}}$.
\begin{table}[t!]
\begin{center}
\caption{\underline{\strut Supervised classification in the VisTex database}}
\label{table:experiment}
\begin{tabular}{cc}
\underline{\strut Classification rule} & \underline{\strut Overall accuracy }\\[0.35cm]
 rule (\ref{eq:mlcrbis}) with  $M = 3$ & $94.3 \pm 0.4 \; \%$ \\[0.13cm] 
 rule (\ref{eq:mlcrbis}) with  $M = 1$  & $86.2 \pm 0.4 \; \%$ \\[0.23cm]
rule (\ref{eq:neighborbis}) with  $M = 3$ & $92.1 \pm 0.5 \; \%$ \\[0.13cm]
rule (\ref{eq:neighborbis}) with  $M = 1$ & $82.5 \pm 0.5 \; \%$ \\[0.23cm]
rule (\ref{eq:wishartclass}) with $M = 3$ & $89.7 \pm 0.8 \; \%$ \\[0.13cm]
rule (\ref{eq:wishartclass}) with $M = 1$ & $84.6 \pm 0.5 \; \%$ \\[0.2cm]
\hline \\[0.2cm]
\end{tabular}
\end{center}
\end{table}

In order to understand Table \ref{table:experiment}, consider the application of the classification rules (\ref{eq:mlcrbis}), (\ref{eq:neighborbis}) and (\ref{eq:wishartclass}), in the context of the present experiment, (the following description is simplified in order to maintain a reasonable length, but a fully detailed version may be found in~\cite{gsi}). 

The experiment reproduces the setting of supervised classification, described in \ref{subsec:laplace}. Indeed, the starting point is a training sequence $\mathcal{T}$, made up of $84 \times 40$ textures. This training sequence is partitioned into disjoint classes, where each class $\mathcal{C}$ consists of textures which belong to the same image, (therefore, these textures represent the same material or object). Since there are $40$ images, $\mathcal{T}$ is partitioned into $40$ classes. 

The textures belonging to each class $\mathcal{C}$ display a diversity of lighting conditions and perspectives. To model this ``in-class diversity," it is assumed that $\mathcal{C}$ may be made up of disjoint clusters $C_1, \ldots, C_{\scriptscriptstyle M\,}$, (eventually, $M = 1$ is considered), where textures belonging to the same cluster display similar lighting conditions and perspectives. 

The training sequence was used to compute the maximum likelihood estimates appearing in the classification rules (\ref{eq:mlcrbis}), (\ref{eq:neighborbis}) and (\ref{eq:wishartclass}). Then, these three classification rules were applied to each one of the $85 \times 40$ classification patches. Simply put, the aim is to retrieve, for each classification patch, the original image from which it was obtained, (and thus the material or object which it represents).

Assume, for simplicity, each texture in the training sequence $\mathcal{T}$ is described by a  parameter $Y \in \mathcal{P}_{m\,}$, (here, $m = 2$ is used). It follows that each class $\mathcal{C}$, within this training sequence, can be identified with a set of data points in $\mathcal{P}_{m\,}$, say $\mathcal{C} = \lbrace Y_1, \ldots, Y_{\scriptscriptstyle N}\rbrace\,$.

 For classification rules (\ref{eq:mlcrbis}) and (\ref{eq:neighborbis}), $Y_1, \ldots,  Y_{\scriptscriptstyle N}$ are modeled as a realisation of a mixture of Gaussian distributions, given by (\ref{eq:mixture}). For classification rule (\ref{eq:wishartclass}), they are modeled as a realisation of a mixture of Wishart distributions, as in~\cite{nielsen}\cite{saintjean}. In either case, the number of mixture components is $M$, the number of clusters contained in $\mathcal{C}$.

The maximum likelihood estimates appearing in rules (\ref{eq:mlcrbis}) and (\ref{eq:neighborbis}) were computed using the EM algorithm of \ref{subsec:em}, and those appearing in rule (\ref{eq:wishartclass}) were computed using the EM algorithm of~\cite{nielsen}\cite{saintjean}, for estimation of mixtures of Wishart distributions. These algorithms were applied to the points $Y_1, \ldots,  Y_{\scriptscriptstyle N\,}$, of each class $\mathcal{C}$, using two choices of the number of mixture components $M$. These are $M = 1$, (which amounts to ignoring in-class diversity), and $M = 3$. It was assumed, (clearly, this is an artificial assumption), that $M$ is the same for all classes contained in the training sequence. 

Table \ref{table:experiment} gives the respective performance of classification rules (\ref{eq:mlcrbis}), (\ref{eq:neighborbis}) and (\ref{eq:wishartclass}). This was obtained by applying these rules to each one of the $85 \times 40$ classification patches. Recall that patches, since they have a suitable scale, can be identified with textures. Therefore, each classification patch was described by a parameter $Y_t \in \mathcal{P}_{m\,}$, which was then used to evaluate the rules (\ref{eq:mlcrbis}), (\ref{eq:neighborbis}) and (\ref{eq:wishartclass}).

Table \ref{table:experiment} shows classification performance in terms of overall-accuracy. As already mentioned, this is the percentage of classification patches which are correctly associated to their original class, (that is, to the image from which they were obtained). Note that rules (\ref{eq:mlcrbis}), (\ref{eq:neighborbis}) and (\ref{eq:wishartclass}) provide additional information, as they associate each classification patch to a specific cluster, a subset of a class. This has no effect on the results shown in the table.

A quick look at Table \ref{table:experiment} shows the following\,: \\[0.1cm]
--- The new classification rule (\ref{eq:mlcrbis}), proposed in \ref{subsec:laplace}, provides significantly better performance than the nearest neighbour rule (\ref{eq:neighborbis}), and the Wishart classifier rule (\ref{eq:wishartclass}), which were proposed by~\cite{congedo} and \cite{wishartclass}, respectively. \\
--- The nearest neighbour rule (\ref{eq:neighborbis}) provides better performance than the Wishart classifier rule (\ref{eq:wishartclass}), when the number of mixture components is $M = 3$. However, this is the other way around, when $M = 1$. \\
--- For each of the three classification rules, performance is improved when training is carried out using $M = 3$, instead of the trivial choice $M = 1$, (which reduces the mixture to a single component). \\[0.1cm]
\indent Validation of the new classification rule (\ref{eq:mlcrbis}), through application to other databases, besides the VisTex database, and also to images acquired from remote sensing projects, is currently ongoing. This will include comparison of this new classification rule to a more comprehensive choice of specialised classification techniques, as may be found in the texture classification literature.
 
\bibliographystyle{IEEEtran}
\bibliography{refs_gausspn}

\section*{Acknowledgements} The authors wish to thank Professor Baba C. Vemuri, of the University of Florida, for his feedback on the paper, and for ongoing collaboration, in view of the future development of related work.

\appendices   

\section{Proof of Proposition \ref{prop:ybar}} \label{appendix0}

This appendix provides the proof of Proposition \ref{prop:ybar}, Paragraph \ref{subsec:gaussianasymp}. A different proof of the first item of the proposition, formula (\ref{eq:gaussbarycentre}), is given in~\cite{vemuri}, (see Theorem 2.3., Page $598$). \\[0.2cm]
\indent --- \textbf{Proof}\,: recall that $\bar{Y}$ is the Riemannian centre of mass of $G(\bar{Y},\sigma)$, if it is a global minimiser of $Y \mapsto \mathcal{E}(Y| \, \bar{Y},\sigma)$. This is the same as saying that $\bar{Y}$ verifies (\ref{eq:gaussbarycentre}). Write $\mathcal{E}(Y)$ in place of $\mathcal{E}(Y| \, \bar{Y},\sigma)$. Proposition \ref{prop:curvature} states that $\bar{Y}$ is the Riemannian centre of mass of $G(\bar{Y},\sigma)$ if $\bar{Y}$ is a stationary point of $\mathcal{E}: \mathcal{P}_m \rightarrow \mathbb{R}_+$.
 
It is now shown this is true. Precisely, denoting $\nabla\mathcal{E}$ the Riemannian gradient of $\mathcal{E}$, it is shown that 
\begin{equation} \label{eq:gradient}
    \nabla\mathcal{E}(\bar{Y}) = 0
\end{equation}
A well-known expression for $\nabla \mathcal{E}$ is the following, see~\cite{afsari},
\begin{equation} \label{eq:gradient1}
  \nabla\mathcal{E}(Y) = -2 \; {\textstyle \bigintsss_{\scriptscriptstyle \mathcal{P}_{m}}} \mathrm{Log}_{\scriptscriptstyle Y} (Z) \, p(Z|\,\bar{Y},\sigma)dv(Z)
\end{equation}
where $\mathrm{Log}_{\scriptscriptstyle Y}$ denotes the Riemannian logarithm mapping, (whose expression is (\ref{eq:log}), given below). 

To prove (\ref{eq:gradient}), note that for all $Y \in \mathcal{P}_{m\,}$,
\begin{equation} \label{eq:one}
{\textstyle \bigintsss_{\scriptscriptstyle \mathcal{P}_{m}}} \; p(Z|\,Y,\sigma)dv(Z) = 1
\end{equation}
since $p(Z|\,Y,\sigma)$, as defined by (\ref{eq:gaussianpdf}), is a probability density. Taking the Riemannian gradient of both sides, it follows
\begin{equation} \label{eq:gradient2}
\nabla_{\scriptscriptstyle Y} \left(\,{\textstyle \bigintsss_{\scriptscriptstyle \mathcal{P}_{m}}} \; p(Z|\,Y,\sigma)dv(Z)\right) = 0
\end{equation}
where $\nabla_{\scriptscriptstyle Y}$ means the gradient is with respect to the variable $Y$. Assume it is possible to carry the Riemannian gradient under the integral. The previous identity (\ref{eq:gradient2}) then becomes 
\begin{equation} \label{eq:gradient3}
 {\textstyle \bigintsss_{\scriptscriptstyle \mathcal{P}_{m}}} \; \nabla_{\scriptscriptstyle Y} \, p(Z|\,Y,\sigma)dv(Z)= 0
\end{equation}
Recall expression (\ref{eq:gaussianpdf}) of $p(Z|\,Y,\sigma)$. Computing the gradient, one has,
\begin{eqnarray} \label{eq:inter}
\nabla_{\scriptscriptstyle Y} \, p(Z|\,Y,\sigma) = \frac{1}{\zeta(\sigma)}\times -\frac{1}{\strut 2\sigma^2}\times 
 \nabla_{\scriptscriptstyle Y} \, d^{\,2}(Z,Y)\times  \exp \left[ - \frac{d^{\,2}(Z,Y)}{2\sigma^2}\,\right] 
\end{eqnarray}
However, the Riemannian gradient of squared Rao's distance is given by, (see~\cite{chavel}, Page $407$), 
\begin{equation} \label{eq:gradlog}
\nabla_{\scriptscriptstyle Y} \, d^{\,2}(Z,Y) = - 2 \, \mathrm{Log}_{\scriptscriptstyle Y} (Z)
\end{equation}
Which means that (\ref{eq:inter}) can be written
$$
\nabla_{\scriptscriptstyle Y} \, p(Z|\,Y,\sigma) = \frac{1}{\strut \sigma^2} \times \mathrm{Log}_{\scriptscriptstyle Y} (Z)\times  p(Z|\, Y,\sigma)
$$
Replacing this expression in (\ref{eq:gradient3}) yields,
\begin{equation} \label{eq:gradient4}
 {\textstyle \bigintsss_{\scriptscriptstyle \mathcal{P}_{m}}} \; \mathrm{Log}_{\scriptscriptstyle Y} (Z)\,  p(Z|\, Y,\sigma) dv(Z)= 0
\end{equation}
which holds for all $Y \in \mathcal{P}_{m\,}$. In particular, putting $Y = \bar{Y}$, and comparing to (\ref{eq:gradient1}), one immediately obtains (\ref{eq:gradient}). This completes the proof of part (i) of the proposition, (the justification for carrying the gradient under the integral in (\ref{eq:gradient2}) can be made using the dominated convergence theorem).

In order to maintain a reasonable length, the proof of part (ii) of the proposition will not be detailed. It uses the same technique as the proof of part (i). Precisely, to obtain (\ref{eq:gaussigma}), differentiate both sides of (\ref{eq:one}) with respect to $\sigma$, and carry the derivative under the integral. It is then possible to conclude by a direct calculation. \hfill $\blacksquare$

\section{Proof of Proposition \ref{prop:normality}} \label{appendix}
This appendix provides the two remaining parts of the proof of Proposition \ref{prop:normality}, Paragraph \ref{subsec:gaussianasymp}. Namely, these two parts are the proofs of Lemmas \ref{lemma:normality1} and \ref{lemma:normality2}. \\[0.2cm]
\indent --- \textbf{Proof of Lemma \ref{lemma:normality1}}\,: to begin, let $\gamma:[0,1] \rightarrow \mathcal{P}_m$ be the geodesic curve connecting $\bar{Y}$ to $\hat{Y}_{\scriptscriptstyle N\,}$. By definition, 
\begin{equation} \label{eq:initv}
  \frac{d}{dt} \, \gamma(0) = \Delta
\end{equation}  
Let $(e_a(t); a = 1, \ldots, p)$, be a parallel orthonormal basis along $\gamma_{\,}$, with $e_a(0) = e_{a\,}$, where $e_a$ is the basis introduced in (\ref{eq:delta}). For the notion of parallel orthonormal basis, see~\cite{chavel}. Then, it is possible to write, 
$$
\nabla \mathcal{E}_{\scriptscriptstyle N}(\gamma(t)) = \sum^p_{a=1} \, \nabla \mathcal{E}^a_{\scriptscriptstyle N}(t) e_a(t)
$$
Indeed, this means that $(\nabla \mathcal{E}^a_{\scriptscriptstyle N}(t); a = 1, \ldots, p)$ are the components of the vector $\nabla \mathcal{E}_{\scriptscriptstyle N}(\gamma(t))$ in the basis $(e_a(t); a = 1, \ldots, p)$. 

From the Taylor development of the functions $\nabla \mathcal{E}^a_{\scriptscriptstyle N}(t)$, 
\begin{equation} \label{eq:taylor1}
\hspace{-0.015cm} \nabla \mathcal{E}^a_{\scriptscriptstyle N}(1) = \nabla \mathcal{E}^a_{\scriptscriptstyle N}(0) + \sum^p_{b=1} \, \nabla^2\mathcal{E}_{\scriptscriptstyle N}(e_b,e_a)(\gamma(0)) \,  \Delta_b + R_{\scriptscriptstyle N} \!\!
\end{equation}
where $R_{\scriptscriptstyle N}$ denotes the remainder. This follows from (\ref{eq:initv}) and from the Taylor development formula in~\cite{chavel} (Page $83$). 

For the left-hand side of (\ref{eq:taylor1}), note that $\gamma(1) = \hat{Y}_{\scriptscriptstyle N}$ is the global minimum of $\mathcal{E}_{\scriptscriptstyle N\,}$. In particular, $\gamma(1)$ is a stationary point of $\mathcal{E}_{\scriptscriptstyle N}$, so $\nabla \mathcal{E}_{\scriptscriptstyle N}(\gamma(1)) = 0$ and
\begin{equation} \label{eq:lhs}
 \nabla \mathcal{E}^a_{\scriptscriptstyle N}(1)  = 0
\end{equation} 
For the right-hand side, note that $\gamma(0) = \bar{Y}$. Therefore, by (\ref{eq:empirical}) and (\ref{eq:gradlog})
\begin{equation} \label{eq:clt1}
  N^{\scriptscriptstyle 1/2} \times \nabla \mathcal{E}_{\scriptscriptstyle N}(\gamma(0)) =  N^{\scriptscriptstyle -1/2}\times -2 \, \sum^{\scriptscriptstyle N}_{n=1} \mathrm{Log}_{\scriptscriptstyle \bar{Y}} Y_n
\end{equation}
Furthermore, using (\ref{eq:gradient}) and (\ref{eq:gradient1}), it follows from the central limit theorem, with $C$ given by (\ref{eq:C}),
\begin{equation} \label{eq:clt2}
   \mathcal{L} \lbrace N^{\scriptscriptstyle 1/2} \times (\nabla \mathcal{E}^1_{\scriptscriptstyle N}(0), \ldots, \nabla \mathcal{E}^p_{\scriptscriptstyle N}(0))  \rbrace \, \Rightarrow \mathcal{N}(0,C)
\end{equation}
In the second term on the right-hand side,
\begin{equation} \label{eq:clt3}
\nabla^2\mathcal{E}_N(e_b,e_a)(\gamma(0)) =   N^{\scriptscriptstyle - 1 } \sum^{\scriptscriptstyle N}_{n=1} \nabla^2 \, d^{\,2}(\bar{Y},Y_n) 
\end{equation}
So, by the law of large numbers, with $H$ given by (\ref{eq:H}),
\begin{equation} \label{eq:lln}
\nabla^2\mathcal{E}_N(e_b,e_a)(\gamma(0)) \longrightarrow H \hspace{0.5cm} \mbox{ as } N \rightarrow \infty
\end{equation}
Finally, the remainder $R_{\scriptscriptstyle N}$ can be written, 
\begin{equation} \label{eq:remainder0}
 R_{\scriptscriptstyle N} = \sum^p_{b=1} \, H^\prime_{ab} \,  \Delta_b 
\end{equation}
where,  by Proposition \ref{prop:consistency}, $H^\prime_{ab}$ verifies the limit
\begin{equation} \label{eq:remainder}
H^\prime_{ab} \longrightarrow 0 \hspace{0.5cm} \mbox{ as } N \rightarrow \infty
\end{equation}
Now, to prove the lemma, it is enough to multiply both sides of (\ref{eq:taylor1}) by $N^{\scriptscriptstyle 1/2}$ and substitute (\ref{eq:lhs}), as well as the limits (\ref{eq:clt2}), (\ref{eq:lln}) and (\ref{eq:remainder}). \hfill $\blacksquare$ \\[0.2cm]
\indent --- \textbf{Proof of Lemma \ref{lemma:normality2}}\,: the aim will be to prove that $H = (1/2\sigma^{\scriptscriptstyle 2})\times C$. Indeed, Lemma \ref{lemma:normality1} states that $\Lambda = H^{-1} \, C \, H^{-1}$. Therefore, $H = (1/2\sigma^{\scriptscriptstyle 2})\times C$ implies $\Lambda = 4\,\sigma^{\scriptscriptstyle 4}\times C^{-1\,}$.

To begin, let $C(v,w)$ be given by 
\begin{equation} \label{eq:Cop}
   C(v,w) \, = \, 4 \times {\textstyle \bigintsss_{\scriptscriptstyle \mathcal{P}_{m}}}  
\left[\, \langle \mathrm{Log}_{\,\scriptscriptstyle \bar{Y}\,}(Z),v\, \rangle_{\scriptscriptstyle \bar{Y}} \, \times \, \langle \mathrm{Log}_{\,\scriptscriptstyle \bar{Y}\,}(Z),w\,\rangle_{\scriptscriptstyle \bar{Y}} \,\right] \,  p(Z|\, \bar{Y},\sigma) dv(Z)
\end{equation} 
for any symmetric $m \times m$ matrices $v$ and $w$. Then, it is clear from (\ref{eq:C}), that $C_{ab} = C(e_a,e_b)$. Recall that $H$ is given by (\ref{eq:H}), in terms of the Riemannian Hessian $\nabla^2 \mathcal{E}(\bar{Y})$. To show that $H = (1/2\sigma^{\scriptscriptstyle 2})\times C$, it will be enough to prove
\begin{equation} \label{eq:hessproof}
  \nabla^2 \mathcal{E}(\bar{Y})(v,w) =  (1/2\sigma^{\scriptscriptstyle 2})\times C(v,w)
\end{equation}
for all symmetric $m \times m$ matrices $v$ and $w$.

In order to obtain this equality, it will be convenient to introduce the following notation. For $Y \in \mathcal{P}_{\scriptscriptstyle m}$ and $\sigma > 0\,$, let 
\begin{equation} \label{eq:lll}
  p(Z|\,Y,\sigma) = e^{\ell(Z|\,Y,\,\sigma)} \hspace{1cm} \ell(Z|\,Y,\sigma) = - \log \zeta(\sigma) - \frac{1}{2\sigma^{\scriptscriptstyle 2}} \, d^{\scriptscriptstyle\, 2}(Z,Y)
\end{equation}
as follows from (\ref{eq:gaussianpdf}). Since $p(Z|\,Y,\sigma)$ is a probability distribution,
$$
{\textstyle \bigintsss_{\scriptscriptstyle \mathcal{P}_{m}}} \, e^{\ell(Z|\,Y,\,\sigma)} \, dv(Z) = 1
$$
for any $Y \in \mathcal{P}_{\scriptscriptstyle m}$. The above integral, as a function of $Y$, being constant, its Riemannian Hessian is equal to $0$. By carrying the Riemannian Hessian under the integral,
\begin{equation} \label{eq:hessproof1}
 {\textstyle \bigintsss_{\scriptscriptstyle \mathcal{P}_{m}}} \, \nabla^2_{\scriptscriptstyle Y}  \, e^{\ell(Z|\,Y,\,\sigma)} \, (v,w) \, dv(Z) = 0
\end{equation}
for all symmetric $m \times m$ matrices $v$ and $w$, where $\nabla^2_{\scriptscriptstyle Y}$ indicates the Riemannian Hessian is with respect to the variable $Y$. The expression under the integral can be evaluated using the following identity, which holds for the Riemannian Hessian of the exponential of any function $\ell$,
\begin{equation} \label{eq:hessproof2}
\nabla^2_{\scriptscriptstyle Y}  \, e^{\ell(Z|\,Y,\,\sigma)} \, (v,w) \, = \, \left \lbrace \,\nabla^2_{\scriptscriptstyle Y} \ell(Z|\,Y,\,\sigma) \, (v,w) \, + \, \langle \nabla_{\scriptscriptstyle Y} \ell(Z|\,Y,\,\sigma),v\rangle_{\scriptscriptstyle Y} \times \langle \nabla_{\scriptscriptstyle Y} \ell(Z|\,Y,\,\sigma),w\rangle_{\scriptscriptstyle Y} \right \rbrace \,  e^{\ell(Z|\,Y,\,\sigma)} 
\end{equation}
where the scalar product notation is that of (\ref{eq:scalarproduct}). 

This identity is shown in the short note following the proof. It is an immediate result of the definition of the Riemannian Hessian~\cite{chavel} (Page $3$ and Page $41$). Using (\ref{eq:lll}), the following expressions are obtained,
\begin{eqnarray} \label{eq:hessproof3}
\label{hessproof31}   \nabla_{\scriptscriptstyle Y}  \, \ell(Z|\,Y,\,\sigma)  = - \frac{1}{2\sigma^{\scriptscriptstyle 2}} \, \nabla_{\scriptscriptstyle Y} d^{\scriptscriptstyle\, 2}(Z,Y)      	= \frac{1}{\sigma^{\scriptscriptstyle 2}} \, \mathrm{Log}_{\scriptscriptstyle Y}(Z) \\[0.15cm]
 \label{hessproof32}  \nabla^2_{\scriptscriptstyle Y}  \, \ell(Z|\,Y,\,\sigma) = - \, \frac{1}{2\sigma^{\scriptscriptstyle 2}} \, \times 
  \nabla^2_{\scriptscriptstyle Y}  \, d^{\,\scriptscriptstyle 2}(Z,Y) \hspace{1.6cm}\,\,
\end{eqnarray}
where the second equality in (\ref{hessproof31}) follows from (\ref{eq:gradlog}). Replacing (\ref{eq:hessproof2}), (\ref{hessproof31}) and (\ref{hessproof32}) into (\ref{eq:hessproof1}) gives
$$
\begin{array}{ll}
- \, \frac{1}{2\sigma^{\mathstrut \scriptscriptstyle 2}} \, \times {\textstyle \bigintsss_{\scriptscriptstyle \mathcal{P}_{m}}} \,  
  \nabla^2_{\scriptscriptstyle Y}  \, d^{\,\scriptscriptstyle 2}(Z,Y)(v,w)\, p(Z|\,Y,\sigma)dv(Z)  & \\[0.2cm]
+ \, \frac{1}{\mathstrut \sigma^{\scriptscriptstyle 4}} \, \times {\textstyle \bigintsss_{\scriptscriptstyle \mathcal{P}_{m}}} \,
\left[\, \langle \mathrm{Log}_{\,\scriptscriptstyle Y\,}(Z),v\, \rangle_{\scriptscriptstyle Y} \, \times \, \langle \mathrm{Log}_{\,\scriptscriptstyle Y\,}(Z),w\,\rangle_{\scriptscriptstyle Y} \,\right]\, p(Z|\,Y,\sigma)dv(Z) & = \; 0
\end{array}
$$
as follows by using the notation of (\ref{eq:lll}). After putting $Y = \bar{Y}$, the last identity is exactly the same as,
\begin{equation} \label{eq:end}
  - \, \frac{1}{2\sigma^{\mathstrut \scriptscriptstyle 2}} \, \times  \nabla^2 \mathcal{E}(\bar{Y})(v,w) \, + \, \frac{1}{\mathstrut 4\sigma^{\scriptscriptstyle 4}} \, \times C(v,w) = 0
\end{equation}
which immediately gives (\ref{eq:hessproof}). This completes the proof of the lemma. \hfill $\blacksquare$ \\[0.2cm]
\indent --- \textbf{Proof of (\ref{eq:hessproof2})}\,: Recall the definition of the Riemannian Hessian of a function $f$, (see~\cite{chavel}, Page $41$),
$$
\nabla^2 f(Y)(v,w) = \langle \nabla_v \, \nabla f(Y),w \,\rangle_{\scriptscriptstyle Y}
$$
where $\nabla_v$ denotes the covariant derivative, (see~\cite{chavel}, Page $3$), and $\nabla f$ is the Riemannian gradient of $f$. If $f(Y) = e^\ell(Y)$, then it follows,
$$
\nabla f(Y) = \nabla \ell(Y) \times e^{\ell(Y)}
$$
Then, by the product formula for the covariant derivative,
$$
\nabla_v \,\nabla f(Y) = \nabla_v \, \nabla \ell(Y) \times e^{\ell(Y)} \, +\,  \nabla \ell(Y) \times \langle \nabla \ell(Y),v\rangle_{\scriptscriptstyle Y} \, e^{\ell(Y)}
$$
forming the scalar product with $w$, this yields
$$
\nabla^2 f(Y)(v,w) = \nabla^2 \ell(Y)(v,w) \times e^{\ell(Y)} \,+\, \langle \nabla \ell(Y),w\rangle_{\scriptscriptstyle Y} \times \langle \nabla \ell(Y),v\rangle_{\scriptscriptstyle Y} \, e^{\ell(Y)} 
$$
which is the same as  (\ref{eq:hessproof2}).

\end{document}